\documentclass[11pt,a4paper,reqno]{amsart}

\usepackage[T1]{fontenc}
\usepackage[utf8]{inputenc}
\usepackage{lmodern}
\usepackage[english]{babel}
\usepackage{microtype}
\usepackage{amsmath,amssymb,amsthm,mathtools}
\usepackage[a4paper,margin=1in]{geometry}
\usepackage[numbers]{natbib}
\usepackage{hyperref}
\input{references-compat.tex}
\usepackage[capitalize,noabbrev]{cleveref}
\usepackage{bbm}
\let\mathbb\mathbbm

\numberwithin{equation}{section}
\newtheorem{teo}{Theorem}
\newtheorem{lem}[teo]{Lemma}
\newtheorem{prop}[teo]{Proposition}
\newtheorem{cor}[teo]{Corollary}
\theoremstyle{definition}
\newtheorem{defi}[teo]{Definition}
\theoremstyle{remark}
\newtheorem{oss}[teo]{Remark}

\title[A probabilistic Weyl law]{A Probabilistic Weyl Law for Quantized Tori with Random Potentials}
\author{Nicola Bertoni}
\address{Department of Mathematics, University of Bologna}
\email{nicola.bertoni4@unibo.it}

\begin{document}
    
    \begin{abstract}
        We study the spectrum of the Weyl quantization of complex-valued functions on the torus subject to small perturbations given by random potentials. We prove that the eigenvalues of the perturbed operator satisfy a Weyl law with high probability in the semiclassical limit, extending the results of T. J. Christiansen and M. Zworski \cite{Probabilistic Weyl Laws for Quantized Tori} and M. Vogel \cite{label_9}.
    \end{abstract}

    \maketitle

    \setcounter{tocdepth}{1} 
    \tableofcontents

    \section{Introduction}
    In this paper we investigate the spectral properties of a certain class of non-selfadjoint pseudodifferential operators. The spectral theory of selfadjoint elliptic operators has a long history, beginning with Weyl's classical result on the asymptotic distribution of the eigenvalues of the Dirichlet Laplacian \cite{Weyl} and extending to broad classes of elliptic differential and pseudodifferential operators. A general formulation is provided by Hörmander's work on the spectral function, which yields Weyl asymptotics for selfadjoint elliptic pseudodifferential operators \cite{HormanderSpectral}; see also Ivrii \cite{Ivrii} for a microlocal treatment of precise spectral asymptotics. In the semiclassical setting, Weyl laws relate the distribution of eigenvalues as \(h\to0\) to phase-space volumes determined by the principal symbol. Fundamental contributions to this theory include the work of Helffer and Robert \cite{HelfferRobert}, while a systematic presentation can be found in the monograph of Dimassi and Sj\"ostrand \cite{label_1}.
 
    By contrast, non-selfadjoint operators do not in general
    satisfy a Weyl law and may exhibit strong spectral instability
    and highly irregular eigenvalue distributions. In particular,
    perturbations of arbitrarily small norm may produce large
    variations in the spectrum. These phenomena have been widely
    studied in the literature; see, for instance, the work of
    Hager and Sj\"ostrand \cite{HagerSjostrand}, who established
    probabilistic Weyl asymptotics for a broad class of
    semiclassical pseudodifferential operators on $\mathbb{R}^d$
    subject to small random perturbations.

   Similar phenomena also arise in discrete models, such as non-selfadjoint Toeplitz matrices. In this setting, works by Davies and Hager \cite{DaviesHager}, Guionnet, Wood, and Zeitouni \cite{GuionnetWoodZeitouni}, and Basak, Paquette, and Zeitouni \cite{BasakPaquetteZeitouni} investigate eigenvalue distributions under deterministic or random perturbations. 

    Moreover, Christiansen and Zworski in \cite{Probabilistic Weyl Laws for Quantized Tori}, as well as Vogel in \cite{label_9}, investigated the corresponding problem for Berezin-Toeplitz quantizations of symbols defined on the $2d$-dimensional torus $\mathbb{T}^{2d}$. In this compact setting, the semiclassical quantization yields operators acting on finite-dimensional Hilbert spaces of dimension $N^d$, where
    $h=(2\pi N)^{-1}$ is the semiclassical parameter which tends to $0$. In these works, the authors introduced perturbations given by complex random full matrices and established probabilistic Weyl-type asymptotics for the perturbed spectra. Oltman \cite{Oltman} extended Vogel's result to Berezin--Toeplitz operators on compact quantizable K\"ahler manifolds.

    Sj\"ostrand studied in \cite{Sjostrand 2} the case of multiplicative random perturbations, proving a probabilistic Weyl law for symbols on the Euclidean phase space $\mathbb{R}^{2d}$ by perturbing their quantizations by small random potentials. Related methods were applied by Becker, Oltman and Vogel \cite{BeckerOltmanVogelGraphene} to random tunneling potentials in the chiral model of twisted bilayer graphene, yielding a probabilistic Weyl law for the associated non-selfadjoint tunneling operator.
    
    Motivated by these results, we investigate the spectral
    behaviour of the torus quantizations considered in \cite{Probabilistic Weyl Laws for Quantized Tori, label_9} under small random potential perturbations, following the strategy introduced by Sj\"ostrand in the Euclidean setting. Although our work is inspired by his ideas in \cite{Sjostrand 2}, their implementation in the toroidal setting requires a discrete analysis, with estimates adapted to a family of matrices whose dimension grows in the semiclassical limit. These adaptations lead to a new result in this setting.
    
    A further distinction with \cite{Probabilistic Weyl Laws for Quantized Tori, label_9} is that the perturbations considered here are of potential type and hence have a diagonal structure in the quantized model. Since such perturbations possess significantly fewer degrees of freedom than full random matrices, the result does not follow directly from Vogel's analysis. In particular, a different singular-value analysis is needed to prove that the smallest singular value of the perturbed matrix is not too small with high probability. This leads to Theorem \ref{teo on quantitative final result}, where we establish a probabilistic Weyl law for this class of random perturbations.
    
    Moreover, in Theorem \ref{empirical measure convergence} we prove the asymptotic convergence of the empirical spectral measures associated to the perturbed operator, similarly to what is done in \cite{Sjostrand and Vogel 2} by Sj\"ostrand and Vogel. This provides a more qualitative counterpart to the preceding quantitative asymptotics.
    
    From a physical perspective, the framework considered here admits a natural interpretation. The Weyl quantization associates to a classical periodic observable $p(x,\xi)$ on the toroidal phase space $\mathbb{T}^{2d}$ a quantum periodic observable acting on a finite-dimensional Hilbert space, i.e. a matrix. When $p$ is complex-valued, the corresponding quantization is non-selfadjoint; thus it is no longer an observable, but rather models the generator of the dynamics of an open quantum system. The perturbations in form of diagonal matrices introduced in this work may be interpreted as random potentials, modelling disorder or impurities in the underlying quantum system. 
    \\
    
\begin{figure}[htbp!] 
        \centering
        \begin{minipage}{0.48\textwidth}
        \centering
        \includegraphics[width=\linewidth]{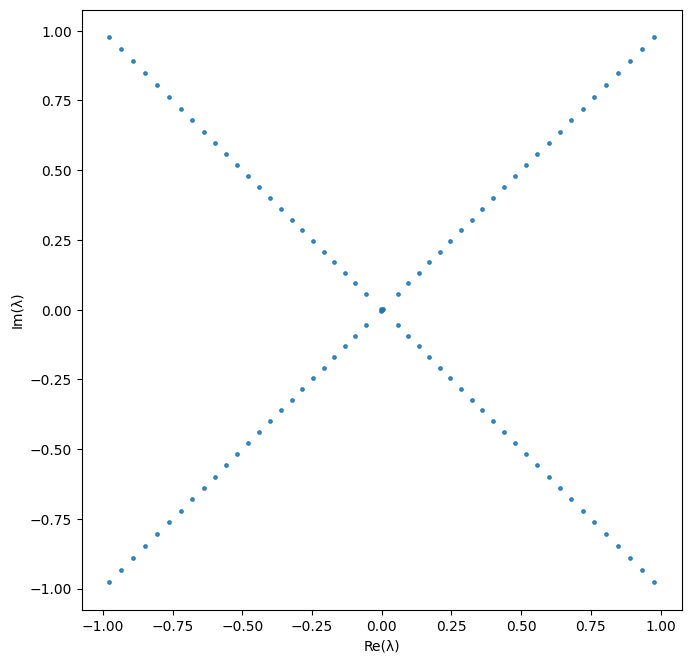}
        \end{minipage}\hfill
        \begin{minipage}{0.48\textwidth}
        \centering
        \includegraphics[width=\linewidth]{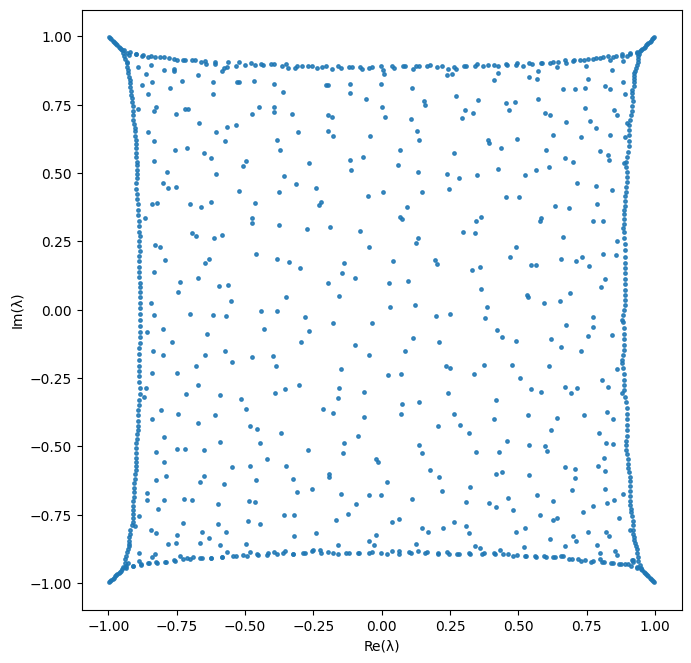}
        \end{minipage}
        \caption{The left hand side shows the spectrum of $p_N$ \eqref{Scott flag op}, while the right hand side shows the spectrum of $p_N+\delta V_q$ for $N=1000$, $V_q$ a complex diagonal Gaussian random matrix and $\delta=N^{-5/2}$.}
        \label{figure simul}
    \end{figure}

    As an illustration of these phenomena, we consider the example of the \textit{Scottish flag operator} given by the symbol
    \begin{equation} \label{symbol scottish flag}
        p(x,\xi) = \cos (2 \pi x) + i \cos (2 \pi \xi), \quad (x,\xi) \in \mathbb{T}^2.
    \end{equation}
    From \eqref{simple quantization}, its quantization is given by
    \begin{equation} \label{Scott flag op}
        p_N = 
        \begin{pmatrix}
            \cos(x_1) & \frac{i}{2} & 0 & 0 & \cdots & \frac{i}{2}\\
            \frac{i}{2} & \cos(x_2) & \frac{i}{2} & 0 & \cdots & 0 \\
            0 & \frac{i}{2} & \cos(x_3) & \frac{i}{2} & \ddots & 0\\
            \vdots & \ddots & \ddots & \ddots & \ddots & \vdots \\
            0 & \cdots & 0 & \frac{i}{2} & \cos(x_{N-1}) & \frac{i}{2} \\
            \frac{i}{2} & 0 & \cdots & 0 & \frac{i}{2} & \cos(x_N) 
        \end{pmatrix}, \quad x_j := \frac{2 \pi j}{N}, \ j = 1, \dots, N.
    \end{equation}
    Figure \ref{figure simul} illustrates two numerical simulations which motivates the approach of random perturbations.
    
    These numerical results are consistent with our analytical results. More precisely, for the Scottish flag operator $p_N$, the following holds
    \begin{teo} \label{Scottish operator theorem}
        Let $Q$ be a diagonal matrix whose diagonal vector is distributed according to the standard Gaussian measure on $\mathbb C^N$ restricted and normalized to the ball $B(0,R) \subset \mathbb{C}^N$, where $ C_1 \leq R \leq C_2 N$, for some constants $C_1, C_2 >0$. Let $\delta = C^{-1}N^{-m-2}$, with $C>0, m \geq \frac{1}{2}$ and let $D$ be an open disc with $\overline{D} \subset (-1,1)^2$. Then, with probability at least $ 1 - \varepsilon_N$,
        \begin{align*}
        \left | \frac{1}{N} \# (\operatorname{Spec}(p_N + \delta Q) \cap D )-  \lambda \, \left({p^{-1}(D)} \right) \right| \leq C N^{-\frac{1}{6}} \quad \mathrm{as} \ N \to \infty,
    \end{align*} 
        where $\varepsilon_N\to0$ as $N\to\infty$.
    \end{teo}
    This theorem is proven in subsection \ref{application section} and it is a simplified version of our main result, which concerns a general class of symbols and random diagonal perturbations. The full statement is given in the following section as Theorem~\ref{teo on quantitative final result}.
    
    \section{Main result} \label{Main result}
    In this paper we will consider semiclassical pseudodifferential operators obtained by the Weyl quantization of smooth periodic symbols on $\mathbb{R}^{2d}$. 
    In particular, this class of operators gives rise to a quantization procedure that maps smooth functions to $N^d \times N^d$ matrices, which are not necessarily selfadjoint:
    \begin{equation} \label{quantization}
        C^\infty(\mathbb{T}^{2d}) \ni p  \mapsto p_N \in \mathcal{L}(\mathbb{C}^{N^d}, \mathbb{C}^{N^d}), \quad N \gg 1.
    \end{equation}
    We will introduce this procedure in more details in Section \ref{Section 2}. As in \cite{Probabilistic Weyl Laws for Quantized Tori}, we observe that when $d=1$, then $\mathbb{T}^2 = S^1_x \times S^1_\xi$ and for functions depending only on one variable it holds
    \begin{align} \label{simple quantization}
        C^\infty(\mathbb{T}^2)\ni f &= f(x) \mapsto f_N = \operatorname{diag}( f(l/N); \ l=0,1,\dots,N-1) \notag \\
        C^\infty(\mathbb{T}^2)\ni g &= g(\xi) \mapsto g_N = \mathcal{F}_N^* \operatorname{diag}(g(l/N); \ l=0,1,\dots,N-1) \mathcal{F}_N,
    \end{align}
    where $\mathcal{F}_N$ is the discrete Fourier transform, defined as in \eqref{Fourier basis}. 

    We are interested in the Weyl quantization of symbols on the torus $\mathbb{T}^{2d}$. We will recall all relevant notions and properties in Section \ref{Section 2}. A function $p \in C^\infty(\mathbb{T}^{2d})$ can be identified with a smooth periodic function on $\mathbb{R}^{2d}$. Hence $p$ is in the symbol class $S(1)$ of smooth functions with a uniform bound in module on all their derivatives.
    We let $h \in ]0,1]$ denote the semiclassical parameter. The $h$-Weyl quantization of a symbol $p \in S(1)$ which may depend on $h$ is given by the linear operator
    $$p^w(x, hD_x; h)u(x) := \frac{1}{(2\pi h)^d} \iint_{\mathbb{R}^{2d}} e^{\frac{i}{h}\langle x-y, \xi \rangle} p \left(\frac{x+y}{2}, \xi;h\right)u(y)  dy d\xi,$$
    acting on $u \in \mathcal{S} (\mathbb{R}^{d})$. The operator $p^w$ is a continuous linear map $\mathcal{S} \rightarrow \mathcal{S}$, $\mathcal{S}' \rightarrow \mathcal{S}'$. 

    We denote by $\mathcal{H}_h^d$ the space of tempered distribution which are $\mathbb{Z}^d$-translation invariant in position and frequency with respect to the semiclassical Fourier transform. The space $\mathcal{H}_h^d$ is $\ne \{ 0\}$ if and only if $h = \frac{1}{2 \pi N},$ for some $\mathbb{N} \ni N >0$, in which case $\dim(\mathcal{H}_h^d) = N^d$, and we can identify $\mathcal{H}_h^d \cong \mathbb{C}^{N^d}$ (see \cite[Lemma 10]{label_9}).

    Taking $p \in C^\infty(\mathbb{T}^{2d})$, we have
    $$p_N := p^w(x,hD_x)_{|\mathcal{H}_h^d}: \mathcal{H}_h^d \rightarrow \mathcal{H}_h^d,$$ 
    which defines the quantization \eqref{quantization}. Now, let $h= \frac{1}{2\pi N}, \mathbb{N} \ni N >0$ and suppose that for $p \in C^\infty(\mathbb{T}^{2d})$ there exist $p_\nu \in C^\infty(\mathbb{T}^{2d})$ so that 
    \begin{equation} \label{expansion}
        p(\rho;h) \sim p_0(\rho) + h p_1(\rho) + \dots \text{ in } S(1),
     \end{equation}
    i.e., $p- \sum_{\nu = 0}^M h^\nu p_\nu \in h^{M+1}S(1)$, for all $ M \in \mathbb{N}$. We will assume also the symmetry condition 
    \begin{equation} \label{symmetry assumption 010}
        p(x,\xi) = p(x,-\xi), 
    \end{equation}
    which implies that
    \begin{equation} \label{symmetry assumption 1}
        p_N^* = \Gamma p_N \Gamma,
    \end{equation}
    where \( \Gamma u := \overline{u} \) denotes the antilinear operator of complex conjugation. Hypothesis \eqref{symmetry assumption 010} plays a fundamental role in the proof of our main result.
    \\
    
    In this article, we study the eigenvalue distribution of
    \begin{equation} \label{object to study}
    p_N + \delta V_q,
    \end{equation}
    where $\delta$ is a small parameter, $V_q = \text{diag}(q_i; i = 1, \dots, N^d)$ and $q \in \mathbb{C}^{N^d}$ is a random vector taking values in the ball $B_{N^d}(0,R) \subset \mathbb{C}^{N^d}$ and distributed according to the probability measure
    \begin{equation}  \label{def of the probability measure 0}
        \mathbb{P}(dq) = C(N) e^{\Phi (q)} \lambda (dq),
    \end{equation}
    where $\lambda$ is Lebesgue measure restricted to this ball and $C(N)$ is the normalization constant. The function $\Phi$ is in $C^1(B_{N^d}(0,R); \mathbb{R})$, it may depend on $N$ and satisfies
    \begin{equation} \label{estimate on phi 0}
        \|\nabla \Phi(q)\|_2 = \mathcal{O}(N^{M_\Phi}),
        \quad \forall\, q \in B_{N^d}(0,R),
    \end{equation}
    for a fixed constant $M_\Phi>0$, uniformly in $q$.
    
    Let $\Gamma \Subset \mathbb{C}$ be an open relatively compact simply connected set with uniform Lipschitz boundary $\partial \Gamma$, according to the definition in Section \ref{section 5.1}. We suppose that there exists $\kappa \in ]0,1]$ such that
    \begin{equation} \label{unif condition}
    V_z(t) := \text{Vol} \{ \rho \in \mathbb{T}^{2d}; |p_0(\rho)-z|^2 \leq t \}= \mathcal{O}(t^\kappa),
    \end{equation}
    uniformly for $z \in \text{neigh}(\partial \Gamma
    )$ and $0 \leq t \ll 1$.
    The first result that we show is the following.

    \begin{teo} \label{teo on quantitative final result}
        Let $p \in C^\infty (\mathbb{T}^{2d})$ satisfy \eqref{expansion} and \eqref{symmetry assumption 010}. Let $\Gamma \Subset \mathbb{C}$ be an open relatively compact simply connected set with uniform Lipschitz boundary $\partial \Gamma$, so that \eqref{unif condition} holds. Let $q \in \mathbb{C}^{N^d}$ be a random vector whose law 
        $\mathbb{P}$ is supported in $B_{N^d}(0,R)$, as defined in \eqref{def of the probability measure 0}, and define $V_q = \mathrm{diag}(q_i ; i = 1, \dots, N^d)$. Let $ e^{-N^{\kappa/2}}\leq \tau_0 \leq(C_\tau N^{-1})^{1/2}$, with $C_\tau >0$ and set
        $$ \varepsilon_0(N) :=  \left(N^{-\frac{d}{2}} + \left( N^{-d} \ln (N) +  N^{-\kappa} \right) \left( \ln \left(\frac{1}{\tau_0}\right) +  \left (\ln (N) \right ) ^2 \right)  \right).$$
        Assume that
        \begin{equation}\label{interval of definition of R}
            C_2N^{\frac{d-\kappa}{2}} \leq R \leq C_3 N^d
        \end{equation}
        where $C>0$ and $C_2>0$ are sufficiently large and independent of $N$. Then, taking $\delta = \frac{\tau_0}{C} N^{-2d}$, $0< r \ll 1$, and $\widetilde{\epsilon} \gg \varepsilon_0 (N)$,
        \begin{align} \label{quantitative final result}
            &\left | \# (\operatorname{Spec}(p_N + \delta V_q) \cap \Gamma) - N^{d} \int_{p_0^{-1}(\Gamma)}  d \rho \right| 
            \leq \mathcal{O}(N^d) \left( 2 \int _{p_0^{-1}(\partial \Gamma + D(0,r))}  d \rho + \mathcal{O}\left(r^{-1} \right) \widetilde{\epsilon} \right),
        \end{align}
        with probability 
        $$\geq 1 - \mathcal{O}(r^{-1}) \varepsilon_0(N)N^{M_\Phi +3d} \exp \left( - \frac{\widetilde{\epsilon}}{C\varepsilon_0(N)}  \right).$$
    \end{teo}

    \begin{oss}
        Condition \eqref{unif condition} is employed throughout the literature in this framework, including in \cite{Probabilistic Weyl Laws for Quantized Tori,HagerSjostrand, label_9}. In particular, it is important for controlling the number of small eigenvalues of $(P-z)^*(P-z)$ for a fixed $z \in \mathbb{C}$. 
        As observed in \cite{Probabilistic Weyl Laws for Quantized Tori},
        if $p$ is real analytic, then \eqref{unif condition} holds for some $\kappa>0$. Similarly, if $p$ is real analytic and $p(\mathbb{T}^{2d}) \subset \mathbb{C}$ has non-empty interior, then for $z \in \mathbb{C}$, the condition $dp_{|p^{-1}(z)} \ne 0$ implies that  \eqref{unif condition} holds with  $\kappa > \frac{1}{2}$.
        Moreover, for $p \in C^\infty (\mathbb{T}^{2d})$ and $z \in \mathbb{C}$,
        \begin{equation} \label{differentials indep}
            dp, d\overline{p} \ \text{ are linearly independent at every point of } \ p^{-1}(z) \Rightarrow \text{\eqref{unif condition} holds with } \, \kappa = 1.
        \end{equation} 
        In particular, we observe that $dp$ and $d \overline{p}$ are linearly independent at a point $\rho$ when 
        $\{ p, \overline{p} \} (\rho) \ne 0$,
        where $\{ a, b \} := \partial_\xi a \partial_x b - \partial_x a \partial_\xi b$ denotes the Poisson brackets on the torus $\mathbb{T}^{2d}.$ The two conditions are equivalent when the dimension $d = 1$, but in dimension $d>1$ this is not true in general. Furthermore, condition \eqref{differentials indep} does not hold when $z \in \partial( p_0(\mathbb{T}^{2d}))$, but it was observed in \cite[Example 12.1]{HagerSjostrand} that for $z \in \mathbb{C}$, if
        $$\forall \ \rho \in p^{-1}(z) : \{ p, \overline{p} \}(\rho) \ne 0 \, \text{ or } \, \{ p, \{ p, \overline{p} \} \}(\rho) \ne 0, \, \text{ then \eqref{unif condition} holds with } \, \kappa = \frac{3}{4}.$$
        
    \end{oss}

    Moreover, following the strategy of \cite{Sjostrand and Vogel 2}, we establish a qualitative result concerning the weak convergence of the empirical spectral measure of the perturbed operator $p_n + \delta V_q$. More precisely, we prove
    \begin{teo} \label{empirical measure convergence}
        Let the symbol $p$, its quantization $p_N$ and the probability measure $\mathbb{P}$ be defined as in the hypotheses of Theorem \ref{teo on quantitative final result}. Taking $0 < \delta \ll 1$ as in Theorem \ref{teo on quantitative final result}, we find
        \begin{equation} \label{weak conv final result}
            \frac{1}{N^d}  \sum _{\zeta \in \sigma(p_N + \delta V_q)} \delta_\zeta\rightharpoonup  (p_0)_*(\lambda(d\rho)), \quad \text{as } \, N \rightarrow + \infty,
        \end{equation}
        in probability, weakly, where $\lambda$ denotes the Lebesgue measure on $\mathbb{T}^{2d}$.
    \end{teo}
    This theorem provides a qualitative counterpart to Theorem~\ref{teo on quantitative final result}. While the latter establishes a probabilistic Weyl law with explicit error estimates, the present result asserts only the asymptotic convergence of the empirical spectral measure, without quantitative bounds.
    
    \subsection{Related results}
    In the case of Toeplitz matrices, numerous results analogous to Theorem \ref{teo on quantitative final result} have been established. Among these, the Jordan matrix has received particular attention. It is a Toeplitz matrix associated with the one-dimensional symbol $p$ defined on $\mathbb{T}$ by
    $ p(\xi)=e^{2\pi i\xi}$. Davies and Hager \cite{DaviesHager} and Sj\"ostrand \cite{Sjostrand 3} showed, under different assumptions on the coupling constant, that with probability close to one, most of the eigenvalues of a randomly perturbed Jordan matrix lie close to the unit circle $S^1$. Sj\"ostrand and Vogel \cite{Sjostr and Vogel 1} subsequently gave a precise description of the eigenvalue distribution in the interior of the disk $D(0,1)$ different from the Weyl law. Using a different approach, \'Sniady \cite{Sniady} established a general convergence result for empirical spectral measures under suitably chosen vanishing Gaussian perturbations. Guionnet, Wood, and Zeitouni \cite{GuionnetWoodZeitouni} later obtained related results for polynomially small Gaussian perturbations. In particular, for the Jordan block matrix, they proved weak convergence in probability of the empirical spectral measure to the uniform probability measure on $S^1$, yielding a result of the same type as Theorem~\ref{empirical measure convergence}. 
    
    Toeplitz matrices associated with more general symbols of the form
    \begin{equation} \label{eq:2.1}
    p(\xi) = \sum_{\nu=-N_-}^{N_+}  p_\nu e^{2\pi i\nu\xi},
    \end{equation}
    where $N_-,N_+\geq 0$, have also been extensively studied. Sj\"ostrand and Vogel in \cite{Sjostrand and Vogel 2} proved a probabilistic Weyl law for randomly perturbed Toeplitz matrices, describing the asymptotic eigenvalue distribution in terms of the symbol $p$. Basak, Paquette and Zeitouni also established in \cite{BasakPaquetteZeitouni} the convergence of the corresponding empirical spectral measure to $p_*\lambda$, where $\lambda$ is the Lebesgue measure. More recently, Bordenave, Capitaine and Chapon \cite{BordenaveCapitaineChapon} provided a detailed analysis of stable and unstable outliers of such randomly perturbed Toeplitz matrices.
    
    More general Toeplitz symbols have also been considered. Sj\"ostrand and Vogel showed in \cite{SjostrandVogelGeneralToeplitz} that for infinite symbols, corresponding to $N_-=N_+=+\infty$ in \eqref{eq:2.1}, and under mild assumptions on $p(S^1)$, an asymptotic Weyl law holds with probability exponentially close to one. They also showed the almost sure convergence of the empirical measure to $p_*\lambda$. No\"el \cite{Noel} extended this type of result to randomly perturbed twisted Toeplitz matrices with low-regularity symbols depending on both position and frequency.
    
    Closer to the setting of the present paper, probabilistic Weyl laws have been obtained for Toeplitz quantizations of symbols. Christiansen and Zworski \cite{Probabilistic Weyl Laws for Quantized Tori} proved that the expected number of eigenvalues in complex domains follows a Weyl law. They conjectured the almost sure convergence of the associated empirical measure in that setting, which was proved by Vogel in \cite{label_9}. In the same paper, Vogel obtained a probabilistic Weyl law and showed that the result of Christiansen and Zworski remained valid for a more general class of random perturbations. More recently, Oltman \cite{Oltman} extended the result of Vogel to the case of Berezin--Toeplitz operators on a K\"ahler manifold.

    In the setting of non-selfadjoint differential and pseudodifferential operators, analogous phenomena were first investigated by Hager \cite{Hager06}, who considered one-dimensional semiclassical differential operators subject to small random perturbations. These results were subsequently extended by Hager and Sj\"ostrand \cite{HagerSjostrand} to a broad class of semiclassical non-selfadjoint pseudodifferential operators in arbitrary dimension, for which they established a Weyl law with high probability under suitable assumptions on the principal symbol and on the random perturbation.

    The case of multiplicative random perturbations was then studied by Sj\"ostrand \cite{Sjostrand 2}, who obtained Weyl asymptotics for non-selfadjoint pseudodifferential operators in arbitrary dimension perturbed by a small random potential. Related spectral instability results for Fredholm operators under random trace-class perturbations were obtained by Becker, Oltman and Vogel \cite{BeckerOltmanVogel}.

    \subsection{Notations}
    In this paper we frequently use the following notations:
    \begin{itemize}
        \item $a \ll b$: $Ca \le b$ for some sufficiently large constant $C > 0$
        \item $f = \mathcal{O}(N)$: there exists a constant $C > 0$ (independent of $N$) such that $|f| \leq CN$.
        \item $f = \mathcal{O}(N^{-\infty})$: for every $M \in \mathbb{N}$, there exists a constant $C_M > 0$, depending on $M$ such that $|f| \leq C_M N^{-M}$.
        \item $\langle \cdot, \cdot \rangle$ is the inner product on any complex vector space under consideration.
    \end{itemize}
    When we want to emphasize that the constant $C > 0$ depends on some parameter $k$, then we write $C_k$, or with the above big-$O$ notation $\mathcal{O}_k(N)$.
    
    \noindent
    Throughout the paper, eigenvalue counts and spectral sums include algebraic multiplicity.

    \subsection*{Acknowledgments}
    The author is very grateful to Martin Vogel for proposing this project and for his guidance and support in the preparation of this paper. The author is especially indebted to Davide Pastorello, who provided constant advice, encouragement, and invaluable insight throughout the development of this paper. The author also wishes to thank the University of Bologna for providing a stimulating and supportive research environment.
    
    The author is a member of the ``Gruppo Nazionale per la Fisica Matematica (GNFM)'' of the ``Istituto Nazionale di Alta Matematica ``Francesco Severi'' (INdAM)''.
    
    \section{Semiclassical calculus and quantization of the torus} \label{Section 2}
    
        In this section, we review some basic notions and results of semiclassical analysis in $ \mathbb{R}^d$, following the exposition in \cite{label_1, label_4, label_12}. We then turn to the study of Toeplitz quantization on the torus, as developed in \cite{Probabilistic Weyl Laws for Quantized Tori} and \cite{Distribution of resonances for open quantum maps}. This framework can be interpreted as a restriction of semiclassical calculus to periodic symbols and to function spaces given by tempered distributions that are periodic both in the spatial variable and in the semiclassical frequency.
        \subsection{Semiclassical quantization}
        We call $m \in C^\infty(\mathbb{R}^{2d}, ]0,+\infty [ )$ an \textit{order function} if there exist constants $C, N_0 > 0$ such that
        $$ m(\rho) \leq C \langle \rho - \eta \rangle^{N_0} m(\eta), \ \ \text{for all} \ \ \rho, \eta \in \mathbb{R}^{2d},$$
        where $ \langle \rho - \eta \rangle := (1+ \|\rho - \eta\|_2^2)^{\frac{1}{2}}$. We define the associated symbol class as  
        \begin{equation} \label{def od S(m)}
            S(m):= \{ p \in C^\infty(\mathbb{R}^{2d}); \ \forall \, \alpha \in \mathbb{N}^{2d}, \ \exists \ C_\alpha >0 \ \ \text{such that} \  |\partial_\rho^\alpha p(\rho)| \leq C_\alpha m(\rho) , \ \forall \, \rho \in \mathbb{R}^{2d} \}.
        \end{equation}
        From now on, we let $h \in ]0,1]$ denote the semiclassical parameter. A symbol $p =p(\rho; h) \in S(m)$ may depend on $h$; in this case we assume that the symbol estimates in the definition \eqref{def od S(m)} hold uniformly in $h$. If a symbol $p \in S(m)$ is of the form $p(\rho; h) = p_0(\rho)+hr(\rho;h), \ r \in S(m)$,
        then $p_0$ is called the \textit{principal symbol} of $p$. Moreover we say that $p$ has the asymptotic expansion
        \begin{equation}
            p \sim p_0 + h p_1 + \dots \text{ in } S(m), \quad p_j \in S(m),
        \end{equation}
        if $p - \sum_{j=0}^N h^j p_j \in h^{N+1}S(m)$, for all $N \in \mathbb{N}$.
        We define the $h$-Weyl quantization of a symbol $p \in S(m)$ as the operator
        \begin{align}
            p^w(x,hD_x) : \mathcal{S}(\mathbb{R}^d) &\rightarrow \mathcal{S}'(\mathbb{R}^d) \notag\\
            u &\mapsto \frac{1}{(2\pi h)^d} \iint_{\mathbb{R}^{2d}} e^{\frac{i}{h}\langle x-y, \xi \rangle}p\left(\frac{x+y}{2}, \xi;h\right)u(y)  dy d\xi, \label{def of Weyl quant}
        \end{align}
        where the $\xi$–integration is understood in the sense of oscillatory integrals. If $p \in S(m)$, then $p^w: \mathcal{S}(\mathbb{R}^{d}) \rightarrow \mathcal{S}(\mathbb{R}^{d})$ and by duality $p^w:  \mathcal{S}'(\mathbb{R}^{d}) \rightarrow \mathcal{S}'(\mathbb{R}^{d})$
        are continuous linear transformations. Moreover, if $p \in S(1)$, then $p^w: L^2(\mathbb{R}^{d}) \rightarrow L^2 (\mathbb{R}^{d})$ is a bounded operator independently of $h$. 

    Given two order functions \(m_1,m_2\) on \(\mathbb{R}^{2d}\), we define the bilinear, continuous \(h\)-Moyal product 
    \begin{equation}\label{def of diesis}
    \#_h:S(m_1) \times S(m_2) \ni (p,q) \mapsto p \#_h q := e^{ihA(D)}(p(x,\xi)q(y,\eta))_{|y=x, \eta = \xi},  \in S(m_1 m_2),
    \end{equation}
    where \(A(D) := \tfrac{1}{2}\,\sigma(D_x,D_\xi; D_y,D_\eta)\), \(\sigma\) denotes the standard symplectic form on \(\mathbb{R}^{2d}\), and we refer to \cite[Theorem 4.8]{label_12} for the integral representation of the quantization of quadratic exponentials. 
    
    From Theorem 4.17 of \cite{label_12}, the $\#_h$-product admits the following asymptotic expansion in the symbol class $S(m_1 m_2)$:
    \begin{equation} \label{semiclassical expansion 01}
    p \#_h q \sim \sum_{k=0}^{+\infty} \frac{1}{k!} \left( ih A(D) \right)^k p(x, \xi) q(y, \eta) \big|_{y=x, \eta=\xi} \ \,   \text{ in } S(m_1 m_2).
    \end{equation}
  
    \begin{prop} [Theorem 4.18 of \cite{label_12}] \label{composition of quantized operators}
        Let $p \in S(m_1), \, q \in S(m_2)$, where $m_1, m_2$ are order functions. Then
        $$p^w(x, hD_x) \circ q^w(x, hD_x) = (p \#_h q)^w(x, hD_x), $$
        as operators mapping $\mathcal{S}(\mathbb{R}^{d}) \to \mathcal{S}(\mathbb{R}^{d})$, and so $\mathcal{S}'(\mathbb{R}^{d}) \to \mathcal{S}'(\mathbb{R}^{d})$ by duality.
    \end{prop}
    A symbol $p \in S(m)$ is called \textit{elliptic} if there exists a constant $ C > 0 $, independent of $h$, such that $|p(x,\xi)| \geq \frac{1}{C} m(x,\xi), \ \forall \ (x,\xi) \in \mathbb{R}^{2d}$.
    
    \subsection{Quantization of the torus}

    We follow the approach of \cite{Probabilistic Weyl Laws for Quantized Tori, Distribution of resonances for open quantum maps}, where the $\mathbb{Z}^d$-periodicity was considered, and of \cite{label_9} for the rescaled torus. 
    Given $\alpha > 0$, we define
    \begin{equation}
    \mathbb{T}^{2d}_\alpha := \mathbb{R}^{2d} \big/ (\alpha^{-1/2} \mathbb{Z})^{2d}.
    \end{equation}
    When $\alpha = 1$, we simply write $\mathbb{T}^{2d} = \mathbb{T}^{2d}_1$. We define the semiclassical Fourier transform by 
    \begin{equation} \label{def of semicl four transform}
        (\mathcal{F}_h u)(\xi) \overset{\text{def}}{=} \frac{1}{(2\pi h)^{d/2}} \int_{\mathbb{R}^{d}} e^{-i \frac{\langle x, \xi \rangle}{h}} u(x) dx, \quad u \in \mathcal{S}(\mathbb{R}^d),
    \end{equation} The operator $\mathcal{F}_h$ maps $\mathcal{S} \to \mathcal{S}$ continuously and can be extended to a continuous function $\mathcal{S}' \to \mathcal{S}'$, mapping $L^2 \to L^2$ unitarily. 
    
    Let $h \ll \alpha \leq 1$ and set $ \widetilde{h} = \frac{h}{\alpha}$. We introduce the space $\mathcal{H}^{d}_{\widetilde{h},\alpha} \subset \mathcal{S}'(\mathbb{R}^d)$, consisting of tempered distributions $u \in \mathcal{S}'(\mathbb{R}^d)$ which are $\alpha^{-1/2} \mathbb{Z}^d$-periodic both in position and in frequency. More explicitly, we require
    \begin{equation} \label{periodicity of u}
    u(x + \alpha^{-1/2} n) = u(x), \quad (\mathcal{F}_{\widetilde{h}} u)(\xi + \alpha^{-1/2} n) = (\mathcal{F}_{\widetilde{h}} u)(\xi), \quad \forall \ n \in \mathbb{Z}^d.
    \end{equation}
    For $\alpha = 1$ we shall simply write $\mathcal{H}^{d}_h := \mathcal{H}^{d}_{h,1}$.  The following lemma characterizes precisely when this space is non-trivial, generalizing the results from the references above.
    \begin{lem} [Lemma 10 of \cite{label_9}]
        Let $h \ll \alpha \leq 1$. Then \begin{equation} \label{basis for H,h,d}
        \mathcal{H}^{d}_{\widetilde{h},\alpha}\ne \{ 0\} \Longleftrightarrow h = \frac{1}{2 \pi N}, \ \text{ for some } N \in \mathbb{N^*} = \mathbb{N} \setminus \{ 0\}.
        \end{equation}
        In this case $\mathcal{H}^{d}_{\widetilde{h},\alpha}$ is a finite-dimensional complex vector space of dimension $N^d$ and
        \begin{equation} \label{basis Q} \mathcal{H}^{d}_{\widetilde{h},\alpha} = \operatorname{span} \left \{ Q_k^\alpha := \frac{1}{(\alpha^{\frac{1}{2}}N)^{\frac{d}{2}}} \sum_{n \in \mathbb{Z}^d} \delta \left(x-\alpha^{-\frac{1}{2}}\left (n+\frac{k}{N} \right) \right); k \in (\mathbb{Z} / N\mathbb{Z})^d  \right \},
        \end{equation}
        where $\delta$ denotes the Dirac distribution.
    \end{lem}
    We observe that the Fourier transform $\mathcal{F}_{\widetilde{h}}$ maps $\mathcal{H}^{d}_{\widetilde{h},\alpha}$ into $\mathcal{H}^{d}_{\widetilde{h},\alpha}$, and can be represented in the basis \eqref{basis Q}, by
    \begin{equation}\label{Fourier basis}
    (\mathcal{F}_{\widetilde{h}})_{n,m} = \frac{e^{-\frac{2\pi i}{N} \langle n, m \rangle}}{N^{d/2}}, \quad n,m \in (\mathbb{Z}/N\mathbb{Z})^{d}.
    \end{equation}
    
    \subsection{Quantization of symbols on the torus and logarithmic estimates} 
    Here, we follow the approach of \cite{label_9}. We assume $h \ll \alpha \leq 1$ and we define an order function on the torus $\mathbb{T}^{2d}_\alpha$ as follows: let $m \in C^\infty( \mathbb{T}^{2d}_\alpha , ]0,+\infty[)$ be such that there exist constants $C_0, N_0 > 0$, independent of $\alpha$, such that
    \begin{align}
    m(\rho) \leq C_0(1 + \|\rho - \mu\|^2_{\mathbb{T}^{2d}_\alpha})^{N_0/2} m(\mu)
    =:  C_0 \langle \rho - \mu \rangle^{N_0}_{\mathbb{T}^{2d}_\alpha} m(\mu), 
    \quad \forall \, \rho, \mu \in \mathbb{T}^{2d}_\alpha, \label{order function on torus}
    \end{align}
    where $\|\rho - \mu\|^2_{\mathbb{T}^{2d}_\alpha} := \inf_{\gamma \in \alpha^{-1/2} \mathbb{Z}^{2d}} \|\rho - \mu + \gamma\|_2^2$. 
    Via the natural projection $\mathbb{R}^{2d} \to \mathbb{T}^{2d}_\alpha$, we may regard $m$ as a $\alpha^{-1/2}\mathbb{Z}^{2d}$–periodic function in $C^\infty(\mathbb{R}^{2d}, \,]0,+\infty[)$. In particular, $m$ is an order function on $\mathbb{R}^{2d}$. 
    
    With this definition, we define the corresponding symbol class
    \begin{equation} \label{def of S(m,alpha)}
    S(m, \alpha) \stackrel{\mathrm{def}}{=}
    \left\{ a \in C^\infty(\mathbb{T}^{2d}_\alpha); \, \forall \, \beta \in \mathbb{N}^{2d} \, \exists \ C_\beta > 0 : 
    |\partial^\beta_\rho a(\rho)| \leq C_\beta m(\rho) , \ \forall \, \rho \in \mathbb{T}^{2d}_\alpha \right\},
    \end{equation}
    where the constants $C_\beta > 0$ are independent of $\alpha$.  
    
    Identifying each symbol in $S(m,\alpha)$ with an $\alpha^{-1/2}\mathbb{Z}^{2d}$–periodic function in $C^{\infty}(\mathbb{R}^{2d})$, we have $S(m,\alpha) \subset S(m)$. This identification will be used throughout, and ensures that the quantization procedure described above applies to all $p \in S(m,\alpha)$. Let $\gamma, \mu \in \alpha^{-1/2}\mathbb{Z}^d$ and consider the unitary operators
    $(\tau_{\gamma} u)(x) := u(x - \gamma), \ (M_\mu u)(x) := e^{\frac{i}{\widetilde{h}} \langle x, \mu \rangle} u(x)$.
    From definition \eqref{def of Weyl quant}, for any $p \in S(m, \alpha)$, $p^w$ commutes with all $\tau_\gamma$ and $M_\mu$. Hence, it follows that
    $$
    p^{w}(x, \widetilde{h} D_{x}) : \mathcal{H}_{\widetilde{h}, \alpha}^{d} \longrightarrow \mathcal{H}_{\widetilde{h}, \alpha}^{d},
    $$
    where $\mathcal{H}_{\widetilde{h}, \alpha}^{d} \subset \mathcal{S}^{\prime}(\mathbb{R}^{d})$. Therefore, we define for $h = \frac{1}{2\pi N}, N \in \mathbb{N}^*$ and $\widetilde{h} := \frac{h}{\alpha}$
    \begin{equation} \label{def of P_N,alpha}
    p_{N, \alpha} \stackrel{\text{def}}{=} p^{w}(x, \widetilde{h} D_{x}) \big|_{\mathcal{H}_{\widetilde{h}, \alpha}^{d}} \in \mathcal{L}(\mathcal{H}_{\widetilde{h}, \alpha}^{d}, \mathcal{H}_{\widetilde{h}, \alpha}^{d}).
    \end{equation}
    When $\alpha = 1$, we simply write $p_{N} := p_{N, 1}$. In particular, we note that $1_{N, \alpha} = \operatorname{Id}_{\mathcal{H}_{\widetilde{h}, \alpha}^{d}}$. 
    
    From the definition \eqref{def of diesis} it follows that if $a,b \in S(m)$ are periodic then $a \#_{\widetilde{h}} b$ is also periodic. Therefore Proposition \ref{composition of quantized operators} applies to symbols $(a, b) \in S(m_{1}, \alpha) \times S(m_{2}, \alpha)$, and yields
    $$
    a_{N, \alpha} \circ b_{N, \alpha} = c_{N, \alpha}, \quad \text{where } c = a \#_{\widetilde{h}} b \in S(m_{1} m_{2}, \alpha).
    $$
    The following lemma characterizes the Hilbert space structure of $\mathcal{H}_{\widetilde{h}, \alpha}^{d}$.
    \begin{lem}[Lemma 11 of \cite{label_9}] \label{inn prod on Hd}
    There exists a unique Hilbert space structure on $\mathcal{H}_{\widetilde{h}, \alpha}^{d}$, determined up to a multiplicative constant, such that every operator $f_{N, \alpha} : \mathcal{H}_{\widetilde{h}, \alpha}^{d} \to \mathcal{H}_{\widetilde{h}, \alpha}^{d}$ with $f \in C^{\infty}(\mathbb{T}_{\alpha}^{2d}; \mathbb{R})$ is self-adjoint. The constant can be chosen so that the basis defined in \eqref{basis Q} is orthonormal. With this choice, the Fourier transform on $\mathcal{H}_{\widetilde{h}, \alpha}^{d}$ defined in \eqref{Fourier basis} is unitary.
    \end{lem}
    From now on, we equip $\mathcal{H}_{\widetilde{h}, \alpha}^{d}$ with the inner product that defines this structure, so that the basis \eqref{basis Q} is orthonormal. Using this basis, we may identify
    \begin{equation} \label{identification of space H}
    \mathcal{H}_{\widetilde{h}, \alpha}^{d} \cong \ell^2((\mathbb{Z}/N\mathbb{Z}) ^d) \cong \mathbb{C}^{N^d}.
    \end{equation}
    From the proof of Lemma \ref{inn prod on Hd} presented in \cite{label_9} we can derive the following result which was presented in the case $\alpha = 1$ in \cite[Lemma 2.5]{Probabilistic Weyl Laws for Quantized Tori}.
    \begin{prop} [Proposition 12 of \cite{label_9}]
        Let $f \in S(m,\alpha)$. Then
        \begin{equation*}
            \operatorname{tr} (f_{N,\alpha}) = (N \alpha)^d \int_{\mathbb{T}_\alpha ^{2d}} f(\rho) d \rho + r_N,
        \end{equation*}
        where for every $k \in \mathbb{N}$, there exists a constant $C_{K,d} > 0$, depending only on $k$ and the dimension $d$, such that
        \begin{equation}
            |r_N| \leq C_{k,d} N^{d-k} \alpha^{d-\frac{k}{2}} \sum_{|\beta| \leq \max (2d+1, k)} \| \partial^\beta f\|_{L^1(\mathbb{T}_\alpha ^{2d})}.
        \end{equation}
    \end{prop}
    We conclude this section with the following proposition.
    \begin{prop}[Proposition 13 of \cite{label_9}] \label{uniform bound for the norm of p_N,alpha}
        Let $N^{-1} \ll \alpha \leq 1$ and $p \in S(1, \alpha)$. Then, there exists a constant $C >0$, independent of $N$ and $\alpha$, such that 
        $$|| p_{N,\alpha} || \leq C.$$ 
    \end{prop}
    These results are presented as in \cite{label_9}, where full proofs can be found.   

    Now, we consider central result concerning log-determinant estimates. The proof is presented in \cite[Proposition 17]{label_9} and it widely makes use of functional calculus for pseudo-differential operators, as presented in \cite[Section 8]{label_1}, and then adapted to symbols in the class $S(m, \alpha)$. 
    
    From now on, we let \( h = (2\pi N)^{-1}, \   N \in \mathbb{N}^* \) and \( N^{-1} \ll \alpha \ll 1 \). Let \( p \in C^{\infty}(\mathbb{T}^{2d}) \) admit an asymptotic expansion \( p \sim p_0 + h p_1 + \ldots \) in \( S(1,1) \). The corresponding Weyl quantization \( P = p^w(x, hD_x; h) \) defines a bounded operator from \( L^2(\mathbb{R}^d) \) to \( L^2(\mathbb{R}^d) \). By \cite[Theorem 4.1]{label_12} $ \left( p^w(x,hD_x) \right)^* = \overline{p}^w(x, hD_x)$. 
    Hence, setting \( Q = P^* P \), by Proposition \ref{composition of quantized operators} there exists a symbol \( q \in S(1,1) \) such that
    \begin{equation}\label{Q is weyl quant of q}
        Q = q^w(x, hD_x; h),
    \end{equation}
    where, from \eqref{semiclassical expansion 01},
    \begin{equation} \label{princ symbol of Q}
    q = \overline{p} \, \#_h \, p \sim q_0 + h q_1 + \ldots \quad \text{in } S(1,1), \quad
    \text{with } q_0 = |p_0|^2.
    \end{equation}
    The following result, established in \cite{label_9}, provides an important log-determinant estimate.
    \begin{prop} [Proposition 17 of \cite{label_9}]  \label{teorema  2.2.2}
        Let $N \in \mathbb{N}^*$ with $N^{-1} \ll \alpha \ll 1$, and let the symbol $q$ be as in \eqref{princ symbol of Q}. Assume that there exists $\kappa \in ]0,1]$ such that the $t$-dependent volume
        \begin{equation} \label{ass for q0}
            V(t):=\lambda \left( \{ \rho \in \mathbb{T}^{2d}; q_0(\rho) \leq t\} \right) = \mathcal{O}(t^{\kappa}), \quad 0 \leq t \ll 1.
        \end{equation}
        Then, for every $\psi \in C^\infty _c(\mathbb{R})$,
        \begin{equation} \label{prima eq big theorem}
            \operatorname{tr}\psi \left( \frac{q_{N}}{\alpha} \right) = N^d \left ( \int \psi \left ( \frac{q_0}{\alpha} \right ) dV(q_0) +\mathcal{O}(N \alpha)^{-1} \alpha^\kappa\right).
        \end{equation}
        Moreover, taking $\chi \in C^\infty _c ([0,+\infty[,[0,+\infty[)$ with $\chi(0) >0$, we have
        \begin{equation} \label{seconda eq big theorem}
            \ln \det \left ( q_N + \alpha\chi \left( \frac{q_N}{\alpha}\right)\right) = N^d \left( \int_{\mathbb{T}^{2d}} \ln q_0(\rho) d\rho + \mathcal{O} \left(\alpha^\kappa \ln \left ( \frac{1}{\alpha}\right )\right) \right) .
        \end{equation}
    \end{prop}
    \begin{oss} \label{last rem chap 1}
        Considering the same assumptions of Proposition \ref{teorema 2.2.2}, the equality \eqref{prima eq big theorem} is particularly useful for estimating the number of eigenvalues of $q_N$ in the interval $[0,\alpha]$. Indeed, let $\psi  \in C_c^\infty (\mathbb{R}, [0,1])$ such that $\psi \equiv 1$ on $[0,1]$ and $\psi \equiv 0$ outside $[  - \frac{1}{2}, \frac{3}{2}]$. Thus, in the right-hand side of \eqref{prima eq big theorem} we obtain
        \begin{align*}
             N^d \left ( \int \psi \left ( \frac{q_0}{\alpha} \right ) dV(q_0) +\mathcal{O}(N \alpha)^{-1} \alpha^\kappa\right) 
             \leq  N^d \left ( \int_0^{\frac{3}{2}\alpha}  dV(q_0) +\mathcal{O}(N \alpha)^{-1} \alpha^\kappa\right)  = \mathcal{O}(N^d \alpha^{\kappa}),
        \end{align*}
        by assumption \eqref{ass for q0}. Furthermore, the left-hand side of \eqref{prima eq big theorem} is bounded from below in the following way
        \begin{align*}
             \operatorname{tr} \left(\psi\left( \frac{q_{N}}{\alpha} \right) \right) \geq  \operatorname{tr} \left(\mathbb{1}_{[0,1]}\left( \frac{q_{N}}{\alpha} \right) \right)
             &= \operatorname{tr} \left(\mathbb{1}_{[0,\alpha]}\left( q_N \right) \right) 
             =  \#\{\zeta \in \sigma (q_{N}): \zeta \in [0,\alpha] \}.
        \end{align*}
        Therefore,
        \begin{equation} \label{corollary of big theorem}
            \#\{\zeta \in \sigma(q_{N}): \zeta \in [0,\alpha] \} = \mathcal{O}(N^d \alpha^{\kappa}).
        \end{equation}
    \end{oss}
    
    \section{Grushin problems} \label{Chapter 2} 
    
    We begin by giving a short overview on Grushin problems. For more general details see for instance \cite{label_8}.
    In a general finite-dimensional setting, we consider an operator 
    $$
    P : \mathcal{H} \longrightarrow \mathcal{H},
    $$
    where \(\mathcal{H}\) is a Hilbert space with $\dim \mathcal{H} < + \infty$. 
    A priori we don't know if $P$ is invertible and we are interested in estimating $\det (P)$. To address this, we introduce the operator
    \begin{equation} \label{grush setting 1}
        \mathcal{P} := 
        \begin{pmatrix}
            P & R_- \\
            R_+  & 0
        \end{pmatrix} 
        : \mathcal{H} \oplus \mathcal{H}_- \longrightarrow \mathcal{H} \oplus \mathcal{H}_+,
    \end{equation}
    where \(\mathcal{H}_+\) and \(\mathcal{H}_-\) are finite-dimensional Hilbert spaces, and 
    $ R_+ : \mathcal{H} \rightarrow \mathcal{H}_+, \  R_- : \mathcal{H}_- \rightarrow \mathcal{H},$
    are suitably chosen operators such that \(\mathcal{P}\) is bijective. 
    If $\dim \mathcal{H}_+ = \dim \mathcal{H}_-$, we may write
    \begin{equation} \label{grush setting 2}
        \mathcal{P}^{-1} = \mathcal{E} := 
        \begin{pmatrix}
            E & E_+ \\
            E_- & E_{-+}
        \end{pmatrix}.
    \end{equation}
    We now turn to the question of the invertibility of \(P\). From the Schur complement formula, it follows that
    \begin{equation} \label{formula cor for P}
        \det (\mathcal{P}) \det(E_{-+}) = \det (P).
    \end{equation}
    Therefore, $P$ is invertible if and only if $E_{-+}$ is invertible. In this case, the corresponding Schur complement is $S = - R_+ P^{-1} R_-$, and moreover
    \begin{equation} \label{expression for E_+}
    E_{-+} = S^{-1} = (- R_+ P^{-1} R_-)^{-1}.
    \end{equation}
    
    \subsection{Grushin problem for the unperturbed operator} \label{Grushin problem for the unperturbed operator}
    In this section we will construct a Grushin problem for the unperturbed operator as in \cite{HagerSjostrand, label_9} making use of its left and right singular vectors. 
    We fix
    $ h = (2 \pi N)^{-1}, \ N \in \mathbb{N}^*$, and $ N^{-1} \ll \alpha \ll 1$, 
    and we consider a symbol \(p \in C^\infty(\mathbb{T}^{2d})\) admitting a semiclassical expansion
    \begin{equation} \label{semiclassical expansion}
    p \sim p_0 + h p_1 + \dots \quad \text{in } \ S(1,1).
    \end{equation}
    We recall the identification $\mathcal{H}_h^d \cong \ell^2((\mathbb{Z}/N\mathbb{Z})^d) \cong \mathbb{C}^{N^d}$ from \eqref{identification of space H},  
    and define $ P := p_N$ 
    as in \eqref{def of P_N,alpha}.
    We fix \(z \in \mathbb{C}\) and consider the operators
    \begin{equation}  \label{def of Q and Q'}
    Q(z) := (P-z)^*(P-z), \quad Q'(z) := (P-z)(P-z)^*.
    \end{equation}
    Let $ 0 \le t_1^2 \le t_2^2 \le \dots \le t_{N^d}^2$
    denote the eigenvalues of \(Q(z)\) with corresponding orthonormal eigenvectors \(e_1, \dots, e_{N^d} \in \mathcal{H}_h^d\). Since \(P\) is a square matrix, we know that
    $ \dim \ker Q(z) = \dim \ker Q'(z) =: N_0$.
    Let \(\{f_1, \dots, f_{N_0}\}\) be an orthonormal basis of \(\ker(P-z)^* = \ker Q'(z)\) and set
    \begin{equation} \label{def of f_i}
    f_i := t_i^{-1} (P-z) e_i, \quad N_0 < i \le N^d.
    \end{equation}
    These vectors are well-defined since \(t_i > 0\) for all \(i > N_0\) and in particular
    $ Q'(z) f_i = t_i^2 f_i.
    $ 
    Hence, \(Q(z)\) and \(Q'(z)\) share the same spectrum \(\{t_i^2\}_{i=1}^{N^d}\). In particular, \(\{f_i\}_{i=1}^{N^d}\) is an orthonormal basis of \(\mathcal{H}_h^d\) of eigenvectors of \(Q'(z)\). Moreover we observe that
    \begin{align}
    (P-z)^* f_i = t_i e_i,   \quad
    (P-z) e_i = t_i f_i, \qquad i = 1, \dots, N^d. \label{properties of the eigenvectors e_j and f_j}
    \end{align} 
    Applying \eqref{Q is weyl quant of q} and \eqref{princ symbol of Q} to $Q(z)$, we obtain
    \begin{equation} \label{princ symb of Q(z)}
        Q(z) = q^w(x, hD_x; h), \qquad q \sim q_0 + h q_1 + \ldots \quad \text{in } S(1,1), \quad
    \text{with } \ q_0 = |p_0-z|^2.
    \end{equation}
    Fix \(M > 0\) so that
    $ 0 \le t_1^2  \leq \dots \leq t_M^2 \leq \alpha < t_{M+1}^2,$
    and let \(\{\delta_i\}_{i=1}^M\) be an orthonormal basis of \(\mathbb{C}^M\). 
    Assuming that the principal symbol \(q_0\) of \(Q(z)\) satisfies \eqref{ass for q0}, Proposition \ref{teorema  2.2.2} and Remark \eqref{last rem chap 1} yield
    \begin{equation} \label{estimate}
    M = \mathcal{O}(N^d \alpha^\kappa).
    \end{equation}
    \begin{oss} \label{rem on uniform on K}
        From the proof of Proposition \ref{teorema 2.2.2} in \cite[Proposition 17]{label_9} it follows that, if \(z\) varies in a compact set \(K \Subset \mathbb{C}\), and the condition \eqref{ass for q0} holds uniformly for \(z \in K\), then the estimate \eqref{estimate} holds uniformly in \(z \in K\).  
    \end{oss}
    We now define the operators
    \begin{align}
        &R_+: \  \mathcal{H}_h^d \longrightarrow \mathbb{C}^M : \quad 
        u \longmapsto \sum_{i=1}^M \langle u, e_i \rangle \, \delta_i, \label{R+} 
        \\
        &R_-: \ \mathbb{C}^M \longrightarrow \mathcal{H}_h^d: \quad
        u_- \longmapsto \sum_{i=1}^M u_-(i) f_i, \label{R-}
    \end{align}
    where \(u_-(i) := \langle u_-, \delta_i \rangle \).  
    The associated Grushin problem is given by
    \begin{equation} \label{def of Pcal(z)}
        \mathcal{P}(z) :=
        \begin{pmatrix}
            P - z & R_- \\
            R_+ & 0
        \end{pmatrix}
        : \mathcal{H}_h^d \oplus \mathbb{C}^M \longrightarrow \mathcal{H}_h^d \oplus \mathbb{C}^M.
    \end{equation}
    We want to show that $\mathcal{P}(z)$ is bijective and derive an explicit formula for its inverse, which we denote by \(\mathcal{E}(z)\). Given \((v, v_+) \in \mathcal{H}_h^d \oplus \mathbb{C}^M\), we wish to solve the equation
    \begin{equation} \label{system}
        \mathcal{P}(z) \begin{pmatrix}
            u \\
            u_-
        \end{pmatrix}
        =
        \begin{pmatrix}
            v \\
            v_+
        \end{pmatrix}.
    \end{equation}
    We write $ u = \sum_{j=1}^{N^d} u_j e_j, \ 
    v = \sum_{j=1}^{N^d} v_j f_j.
    $ 
    Similarly, we write \(u_-\) and \(v_+\) in the basis \(\{\delta_i\}_{i=1}^M\), with coefficients \(\{u_-(i)\}_{i=1}^M\) and \(\{v_+(i)\}_{i=1}^M\) respectively.
    Substituting into \eqref{system}, we obtain
    \begin{equation*}
    \begin{cases}
    \sum_{i=1}^{N^d} t_i u_i f_i + \sum_{j=1}^{M} u_-(j) f_j = \sum_{k=1}^{N^d} v_k f_k, \\
    \sum_{i=1}^M u_i \delta_i = \sum_{j=1}^M v_+(j) \delta_j.
    \end{cases}
    \end{equation*}
    By linear independence of \(\{f_j\}_{j=1}^{N^d}\) and \(\{\delta_j\}_{j=1}^M\), this system is equivalent to
    \begin{equation} \label{system rewritt}
    \begin{cases}
    t_j u_j = v_j, & j = M+1, \dots, N^d, \\[0.5em]
    \begin{pmatrix}
        t_i & 1 \\
        1 & 0
    \end{pmatrix}
    \begin{pmatrix}
        u_i \\
        u_-(i)
    \end{pmatrix}
    =
    \begin{pmatrix}
        v_i \\
        v_+(i)
    \end{pmatrix}, & i = 1, \dots, M.
    \end{cases}
    \end{equation}
    Since
    \begin{equation} \label{inv matr}
    \begin{pmatrix}
        t_i & 1 \\
        1 & 0
    \end{pmatrix}^{-1}
    =
    \begin{pmatrix}
        0 & 1 \\
        1 & -t_i
    \end{pmatrix}, \qquad i=1, \dots, M, 
    \end{equation}
    and $t_j >0$, for all $j \geq M+1$, we deduce that $\mathcal{P}(z)$ is bijective with inverse
    \begin{equation} \label{equ inv P}
        \mathcal{P}^{-1}(z) = \mathcal{E}(z) :=
        \begin{pmatrix}
            E(z) & E_+(z) \\
            E_-(z) & E_{-+}(z)
        \end{pmatrix},
    \end{equation}
    where
    \begin{align}
        &E(z) = \sum_{i=M+1}^{N^d} \frac{1}{t_i} \, e_i f_i^*,  &
        &E_+(z) = \sum_{i=1}^M e_i \, \delta_i^*, \label{definition of E, E+} \\
        &E_-(z) = \sum_{i=1}^M \delta_i \, f_i^*,  &
        &E_{-+}(z) = -\sum_{i=1}^M t_i \, \delta_i \delta_i^*. \label{definition of E-, E-+}
    \end{align}
    Moreover,
    \begin{equation} \label{diseq matr}
        \| E(z) \| \leq \frac{1}{t_{M+1}}, \qquad
        \| E_{\pm}(z) \| = 1, 
        \qquad 
        \|E_{-+}(z)\| \leq t_{M}.
    \end{equation}
    It follows from \eqref{system rewritt} that
    \begin{equation} 
        |\det  (\mathcal{P}(z)) |^2 
        = \prod_{i = M+1}^{N^d} t_i^2 
        = \alpha^{-M} \prod_{i = 1}^{N^d} \mathbf{1}_\alpha (t_i^2), 
    \end{equation}
    where $\mathbf{1}_\alpha(x) := \max(x, \alpha)$.  
    Hence,
    \begin{equation} \label{eq det}
        |\det (\mathcal{P}(z))|^2 = \alpha^{-M} \det ( \mathbf{1}_\alpha(Q(z))).
    \end{equation}
    We can apply Schur complement formula, as in \eqref{formula cor for P}, obtaining 
    $$
    \ln |\det(P-z)| = \ln |\det (\mathcal{P}(z))| + \ln |\det (E_{-+}(z))|.
    $$
    Now let $\chi \in C_c^\infty([0,\infty))$ be a smooth cutoff function supported in $[0,2]$, with $0 \leq \chi \leq 1$ and $\chi \equiv 1$ on $[0,1]$. For every $x \geq 0$, we have the inequalities
    \begin{equation} \label{ineq cutoff}
        x + \frac{\alpha}{4} \chi \left( \frac{4x}{\alpha} \right) \leq \mathbf{1}_\alpha(x) \leq x + \alpha \chi \left( \frac{x}{\alpha} \right).
    \end{equation}
    As shown in \eqref{princ symb of Q(z)}, the principal symbol of the operator $Q(z)$ is given by $|p_0 - z|^2$. By combining Proposition \ref{teorema  2.2.2} with the estimates \eqref{eq det}, \eqref{ineq cutoff}, and \eqref{estimate}, we obtain 
    \begin{align}  \label{formula per P cal}
        \ln |\det \mathcal{P}(z)|^2 
        = \ln \big| \det (\mathbf{1}_\alpha(Q(z))) \big| + M \ln  \frac{1}{\alpha}
        = N^d \left( \int_{\mathbb{T}^{2d}} \ln |p_0(\rho) - z|^2 \, d\rho + \mathcal{O}\!\left( \alpha^\kappa \ln  \frac{1}{\alpha} \right) \right).
    \end{align}
    By Remark \ref{rem on uniform on K}, \eqref{formula per P cal}
    holds uniformly for $z\in K$, for every compact set
    $K\Subset\mathbb C$ on which the volume condition
    \eqref{ass for q0} holds uniformly.
    \begin{oss}
        If $M =0$, then $t_1^2 > \alpha$, and therefore all the singular values of $P-z$ are strictly positive. In particular, $P-z$ is bijective. Hence, there is no need to introduce the auxiliary spaces and the operators $R_+$ and $R_-$, and the Grushin problem reduces to 
        $$\mathcal{P}(z) = P-z.$$
        Therefore,
        \begin{align} 
            \ln |\det \mathcal{P}(z)|^2 =  \ln |\det (P-z)|^2
            = \ln \big| \det (Q(z)) \big|.
        \end{align}
        It follows that \eqref{formula per P cal} remains valid in the case $M=0$.
    \end{oss}

    \subsection{Grushin problem for the perturbed operator} \label{Grushin for perturbed op}
    We consider a deterministic linear perturbation $Q : \mathcal{H}_h^d \rightarrow \mathcal{H}_h^d$ (i.e. an $N^d \times N^d$ matrix since $\mathcal{H}_h^d \simeq \mathbb{C}^{N^d}$ from \eqref{identification of space H}). We study the properties of the operator $P=p_N$ introduced in Section \ref{Grushin problem for the unperturbed operator}, under the perturbation $Q$, using the associated Grushin problem. Our goal is to analyze the spectrum of
    $P^\delta := P + \delta Q,$ with $ \ 0 \leq \delta \ll 1$. For this purpose we set up a Grushin problem for the perturbed operator. 
    
    Fix $z \in \mathbb{C}$. Let the operators $R_{\pm}$ be as in \eqref{R+} and \eqref{R-} and set
    \begin{equation} \label{primo grushin perturbed}
        \mathcal{P}^\delta (z) := \begin{pmatrix}
            P^\delta - z & R_- \\
            R_+ & 0
        \end{pmatrix} :
        \mathcal{H}_h^d \oplus \mathbb{C}^M \longrightarrow \mathcal{H}_h^d \oplus \mathbb{C}^M,
    \end{equation}
    We suppose that $\delta \geq 0$ is such that
    \begin{equation} \label{diseq Q}
    \delta \, \alpha^{-1/2} \, \|Q \|  \leq \frac{1}{2}.
    \end{equation} 
    Let $\mathcal{E}(z) = \mathcal{P}^{-1}(z)$, as defined in \eqref{def of Pcal(z)} and \eqref{equ inv P}. We compute
    \begin{equation*}
    \mathcal{P}^\delta(z) \, \mathcal{E}(z) = I + \delta
    \begin{pmatrix}
     Q E(z) &  Q E_+(z) \\
    0 & 0
    \end{pmatrix}.
    \end{equation*}
    Applying \eqref{diseq matr} and \eqref{diseq Q}, we get
    \begin{align*}
      \left \|  - \delta   \begin{pmatrix}
     Q E(z) &  Q E_+(z) \\
    0 & 0
    \end{pmatrix}\right \|  \leq \delta \| Q \| (\| E(z) \| + \| E_+(z) \|) 
    \leq  \delta \alpha^{-1/2} \| Q \| + \delta \| Q \| < 1.
    \end{align*}
    Therefore, by the Neumann series, we know that $\mathcal{P}^\delta(z) \, \mathcal{E}(z)$
    is bijective with inverse
    \begin{align*}
    (\mathcal{P}^\delta(z) \mathcal{E}(z))^{-1} =   I + 
    \sum_{k=1}^{+\infty} (-\delta)^k  
    \begin{pmatrix}
    (Q E(z))^k & (Q E(z))^{k-1} Q E_+(z) \\
    0 & 0
    \end{pmatrix},
    \end{align*}
    whose norm is $\leq \mathcal{O}(1)$. Consequently, $\mathcal{P}^\delta(z)$ is bijective with inverse
    \begin{align}
        \mathcal{E}^\delta(z) = \mathcal{E}(z) + \sum_{k=1}^{+\infty} (-1)^k
        \begin{pmatrix}
        E(\delta Q E)^k & (\delta E Q )^k E_+\\
        E_-(\delta Q E)^k & E_-(\delta Q E)^{k-1} (\delta Q) E_+
        \end{pmatrix} =: 
        \begin{pmatrix}
        E^\delta(z) & E_+^\delta(z) \\
        E_-^\delta(z)  & E_{-+}^\delta(z) 
        \end{pmatrix}. \label{Sviluppo di Neumann di cal(E)}
    \end{align}
    Hence, by \eqref{diseq matr} and \eqref{diseq Q}, we obtain the following inequalities
    \begin{align}
        &\| E^\delta(z) \| = \left \| E \left( I + \sum_{k=1}^{+\infty} (-1)^k
        (\delta Q E)^k \right) \right \| 
         \leq  \frac{2}{t_{M+1}}. \label{estimate for E^delta} \\
        &\| E_\pm^\delta(z) \|  \leq 2; \label{estimate for E+-^delta} \\
        &\| E_{-+}^\delta(z) - E_{-+}\| 
        \leq   2 \delta \| Q\| \leq \alpha^{1/2}  \label{estimate for E-+^delta}.
    \end{align}
    Applying \eqref{formula cor for P} to $\mathcal{P}^\delta(z) $ and $\mathcal{E}^\delta(z)$, we obtain
    \begin{equation} \label{important sum for the final estimate}
        \ln | \det(P^\delta-z)| = \ln | \det (\mathcal{P}^\delta(z)) |+ \ln | \det (E_{-+}^\delta(z))|. 
    \end{equation}
    By Jacobi's formula 
    $$\frac{d}{d\delta} \det( \mathcal{P}^\delta(z)) = \operatorname{tr} \left(\operatorname{adj}(\mathcal{P}^\delta(z)) \frac{d}{d \delta} \mathcal{P}^\delta(z) \right),$$
    where $\operatorname{adj}(\mathcal{P}^\delta(z))$ denotes the adjugate of $\mathcal{P}^\delta(z)$. Therefore, by \eqref{primo grushin perturbed}, \eqref{Sviluppo di Neumann di cal(E)},
    \begin{align}
        \frac{d}{d\delta} (\ln |\det (\mathcal{P}^\delta(z))| ) &= \Re\left( \frac{1}{\det (\mathcal{P}^\delta(z))} \operatorname{tr} \left(\operatorname{adj}(\mathcal{P}^\delta(z)) \frac{d}{d \delta} \mathcal{P}^\delta(z) \right) \right) = \Re (\operatorname{tr} (E^\delta(z) Q)). \label{formula with d/ddetla}
    \end{align}
    Applying \eqref{estimate for E^delta} and \eqref{formula with d/ddetla}, we obtain
    \begin{align}
        \left | \ln | \det (\mathcal{P}^\delta(z)) | - \ln | \det (\mathcal{P}(z)) | \right | = \left |  \int _0^\delta \frac{d}{d\tau} \left(\ln | \det( \mathcal{P}^\tau(z))| \right) d\tau \right |  \leq \mathcal{O}(\delta \alpha^{-1/2} \| Q \|_{\operatorname{tr}}). \label{estimate per P cal}
    \end{align} 
    Hence, combining \eqref{formula per P cal} and \eqref{estimate per P cal}, we have the following estimate
    \begin{align}
         \ln | \det (\mathcal{P}^\delta (z))| \leq N^d \left( \int_{\mathbb{T}^{2d}} \ln |p_0(\rho) - z| d\rho + \mathcal{O}\left( \alpha^\kappa \ln \left( \frac{1}{\alpha} \right)\right) +  \mathcal{O}(\delta N^{-d}\alpha^{-1/2} \| Q \|_{\operatorname{tr}})\right). \label{estimate for mathcal(P)delta}
    \end{align}
    Moreover, under the assumption \eqref{diseq Q}, \eqref{estimate for E-+^delta} yields 
    $ \|E_{-+}^{\delta}(z)\| \leq \alpha^{1/2} + t_M =  \mathcal{O}(\alpha^{1/2}),
    $ 
    which, in view of \eqref{estimate}, implies
    \begin{equation} \label{estimate for E-+ delta}
    \ln |\det E_{-+}^{\delta}(z)| \leq M \ln (\| E_{-+}^{\delta}(z)\| )\leq \mathcal{O}(N^d \alpha^\kappa) | \ln \alpha |.
    \end{equation}
    
    We conclude this section with a general estimate on the singular values arising in Grushin problems. Throughout the remainder of the article, for a matrix $A \in \mathbb{C}^{M \times M}$ with singular values ordered increasingly as $t_1(A) \leq \dots \leq t_M(A),$ we use the notation
    $$s_j(A) := t_{M - j + 1}(A), \qquad j=1,\dots,M.$$
    \begin{lem} \label{estimates for sing values}
        Let $\mathcal{H}$ be an $N$-dimensional complex Hilbert space, and let $N \geq M > 0$. Suppose that
        \[
        \mathcal{P} = \begin{pmatrix}
        P & R_- \\
        R_+ & 0
        \end{pmatrix} : \mathcal{H} \oplus \mathbb{C}^M \to \mathcal{H} \oplus \mathbb{C}^M
        \]
        is bijective with inverse
        \[
        \mathcal{E} = \begin{pmatrix}
        E & E_+ \\
        E_- & E_{-+}
        \end{pmatrix}.
        \]
        Let $0 \leq t_1(P) \leq \ldots \leq t_N(P)$ denote the eigenvalues of $(P^* P)^{1/2}$, and let $0 \leq t_1(E_{-+}) \leq \ldots \leq t_M(E_{-+})$ denote the eigenvalues of $(E_{-+}^* E_{-+})^{1/2}$. Then, for $m = 1, \dots, M$, we have
        \begin{equation} \label{aux result for sing val}
        \frac{t_m(E_{-+})}{\|E\| t_m(E_{-+}) + \|E_-\| \|E_+\|} \leq t_m(P) \leq \|R_+\| \|R_-\|  t_m(E_{-+}).
        \end{equation}
    \end{lem}

    \begin{proof}
        We have already seen from \eqref{formula cor for P} that $P$ is invertible if and only if $E_{-+}$ is invertible and in that case, by \eqref{expression for E_+}
        \begin{equation} \label{results for P and E}
            P^{-1} = E - E_+ E_{-+}^{-1}E_-, \qquad E_{-+}^{-1} = - R_+P^{-1}R_-.
        \end{equation}
       We adopt a decreasing notation for the singular values. In particular, we denote the singular values of $P$ as $0 \leq s_N(P) \leq \ldots \leq s_1(P)$ and the singular values of $E_{-+}$ as $0 \leq s_M(E_{-+}) \leq \ldots \leq s_1(E_{-+})$. Thus, 
        if we assume $P$ to be invertible, we have
        \begin{align} 
            &s_n(P^{-1}) = \frac{1}{t_n(P)}, \quad n = 1, \dots, N; \label{invers sing val of P} \\ 
            &s_m(E_{-+}^{-1}) = \frac{1}{t_m(E_{-+})}, \quad m = 1, \dots, M. \label{inverse sing val of E-+}
        \end{align}
        Recall also that $s_1(A) = \| A\|$ for every matrix $A $. 
        Applying Ky Fan’s inequalities (\cite[Corollary 2.2]{Gohberg}) together with the identities \eqref{results for P and E}, it follows that
        \begin{align}
            s_m(P^{-1}) = s_m(E - E_+E_{-+}^{-1}E_- ) 
            \leq \| E \| + \| E_+\| \| E_-\|s_m(E_{-+}^{-1}), \quad  m =1, \dots, M.
        \end{align}
        Therefore, invoking \eqref{invers sing val of P} and \eqref{inverse sing val of E-+} yields
        \begin{align*}
            t_m(P) \geq \frac{1}{ \| E \| + \| E_+\| \| E_-\| \frac{1}{t_m(E_{-+})}} 
            = \frac{t_m(E_{-+})}{ \| E \| t_m(E_{-+}) + \| E_+\| \| E_-\| }, \quad  m =1, \dots, M,
        \end{align*}
        which is the desired lower bound. 
        Moreover, from the second identity in \eqref{results for P and E}, we get
        $$
        s_m(E_{-+}^{-1})=s_m(R_+P^{-1}R_-)
        \ \le\ s_1(R_+)\,s_m(P^{-1})\,s_1(R_-)
        = \|R_+\|\,\|R_-\|\,s_m(P^{-1}).
        $$
        Taking the reciprocals and using \eqref{invers sing val of P} and \eqref{inverse sing val of E-+},
        $$
        t_m(P)\ \le\ \|R_+\|\,\|R_-\|\,t_m(E_{-+}),
        \qquad m=1,\dots,M.
        $$
        Assume now that $P$ is not invertible. We introduce the perturbed operator $P^\varepsilon := P + \varepsilon X$, where $\|X\| \leq 1$ and $0 < \varepsilon \ll 1$, so that $P^\varepsilon$ is bijective. 
        By a Neumann series argument, the associated Grushin problem $\mathcal{P}^\varepsilon$ remains invertible, with inverse $\mathcal{E}^\varepsilon$. Consequently, we may apply \eqref{aux result for sing val} to this perturbed problem. Moreover it can be shown that $ \| P^\varepsilon - P\| \rightarrow 0, \  \| \mathcal{E}^\varepsilon - \mathcal{E}\| \rightarrow 0$, as $ \varepsilon \rightarrow 0$. Hence it follows that \eqref{aux result for sing val} also holds in the case where $P$ fails to be invertible. 
    \end{proof}
          
    \section{Estimates for the perturbed operator} \label{Chap 3}
    In this section we construct a suitable perturbation for the operator $P = p_N$ defined in  Section \ref{Grushin problem for the unperturbed operator}, in order to obtain quantitative lower bounds for the small singular values of the perturbed operator. Our approach adapts the method used by Sj\"ostrand in \cite{Sjostrand 2} to the toroidal phase space $\mathbb{T}^{2d}$.
    We proceed in two stages. First, we build a preliminary perturbation that already furnishes effective lower bounds for a larger set of the singular values of the perturbed operator. Second, following Sj\"ostrand’s iterative scheme, we propagate the estimate down to the bottom of the singular-value scale, thereby controlling even the lowest singular value of the perturbed operator. 
    
    \subsection{Construction of a potential} \label{Section from Sjostrand}
    We fix $N \in \mathbb{N}^*$ and consider the following results.
    \begin{prop}\label{Prop 01}
        Let $M \in \mathbb{N}, \ 1 \leq M\leq N^d$ and let $e_1, \dots, e_M \in \mathbb{C}^{N^d}$ be linearly independent vectors, such that $e_i = ( e_i(n) )_{n=1}^{N^d}$, for each $i = 1, \dots, M$. Then there exist $M$ different points $n_1, \dots, n_M \in \{ 1, \dots, N^d\}$ such that the vectors $\vec{e}(n_1), \dots, \vec{e}(n_M)$ are linearly independent in $\mathbb{C}^M$, where
        \begin{equation} \label{notation for vector e}
            \vec{e}(n) := \begin{pmatrix}
                e_1(n) \\
                e_2(n) \\
                \vdots \\
                e_M(n)
            \end{pmatrix}, \quad \forall \ n \in \{ 1, \dots, N^d\}.
        \end{equation}
    \end{prop}
    \begin{proof}
        We define $E:= \operatorname{Span}\{\vec{e}(n); \ n =  1, \dots, N^d \} \subset\mathbb{C}^M$ and we suppose that $E \ne \mathbb{C}^M$. Hence, there exists $\boldsymbol{0} \ne \lambda = (\lambda_1, \dots, \lambda_M) \in \mathbb{C}^M$ such that $\lambda \in E^\perp$, i.e.
        $ 0 = \langle \lambda , \vec{e}(n) \rangle$, for all $ n \in \{ 1, \dots, N^d\}$,
        and consequently
        $  \boldsymbol{0} =   \sum_{j=1}^M \overline{\lambda_j} e_j. $
        This implies that $e_1, \dots, e_M$ are linearly dependent vectors, which is in contradiction with the hypothesis. Thus, $E = \mathbb{C}^M$.
    \end{proof}

    \begin{prop} \label{Prop 02}
        Let $M \in \mathbb{N}, \ 1 \leq M \leq N^d$ and let $\{e_1, \dots, e_M\}$, $\{f_1, \dots, f_M\}$ be two linearly independent families in $\mathbb{C}^{N^d}$. For both of them, we adopt the notations introduced in Proposition \ref{Prop 01}. We assume that we can find $M$ points $n_1, \dots, n_M \in \{ 1, \dots, N^d\}$ such that both $\left\{ \vec{e}(n_i) \right\}_{i=1}^M$ and $\left \{ \vec{f}(n_i) \right\}_{i=1}^M$ are linearly independent families in $\mathbb{C}^M$. We define
        \begin{align} 
            B: \mathbb{C}^M  \longrightarrow \mathbb{C}^M:  \quad u  \longmapsto \sum_{j = 1}^M  \left\langle u ,\vec{f}(n_j) \right \rangle \vec{e}(n_j). \label{def of M}
        \end{align}
        Then $B$ is bijective.
    \end{prop}
    \begin{proof}
        Let $u \in \operatorname{Ker} \, B$. Since $\vec{e}(n_1), \dots, \vec{e}(n_M)$ form a basis in $\mathbb{C}^M$, we have $ \langle u, \vec{f}(n_j) \rangle = 0$, for all $ j = 1, \dots, M. $ 
        Since $\vec{f}(n_1), \dots, \vec{f}(n_M)$ also form a basis in $\mathbb{C}^M$, it follows that $u = \boldsymbol{0}$. Hence, $B$ is injective and therefore bijective. 
    \end{proof}
    We observe that, under the same assumptions of Proposition \ref{Prop 02}, the matrix associated with $B$, expressed in the canonical basis, is given by
    \begin{equation} \label{entries of M}
    (B)_{j,k} = \sum _{\nu = 1}^M e_j(n_\nu) \overline{f_k (n_\nu)}, \quad j,k  = 1. \dots, M.
    \end{equation}
    Hence, defining
    \begin{align}
        E:= (\vec{e}(n_1), \dots, \vec{e}(n_M)) \in \mathbb{C}^{M \times M}, \qquad F:= (\vec{f}(n_1), \dots, \vec{f}(n_M)) \in \mathbb{C}^{M \times M} \label{def of E and F},
    \end{align}
    we have $ B = EF^*$.
    \begin{lem} \label{lem sjostr}
        Let $M \in \mathbb{N}, \ 1 \leq M\leq N^d$ and let $e_1, \dots, e_M$ be as in Proposition \ref{Prop 01}. Let $L \subset \mathbb{C}^M$ be a linear subspace of dimension $S-1$, for some $1 \leq S \leq M$. Then, there exists $n \in \{ 1, \dots, N^d\}$ such that
        \begin{equation}\label{usef diseq lem}
            \operatorname{dist} \, \left(\vec{e}(n), L \right)^2 \geq \frac{1}{N^d} \operatorname{tr} \,((I - \pi_L) \, \mathcal{E}_{N^d}), 
        \end{equation}
        where $\mathcal{E}_{N^d} := \left( \langle e_j , e_k\rangle \right)_{1 \leq j,k \leq M} $ and $\pi_L$ is the orthogonal projection from $\mathbb{C}^M$ onto $L$. 
    \end{lem}
    \begin{proof}
    Let $\{\nu_1, \dots, \nu_M\}$ be an orthonormal basis of $\mathbb{C}^M$ chosen so that $L = \operatorname{span} \, \{\nu_1, \dots, \nu_{S-1} \}$ ($L = \{ 0\}$ when $S=1$). We write $ \nu_l = (
            \nu_{i,l})_{i =1,\dots,M} \in \mathbb{C}^M,$ for $l = 1, \dots, M$.
    By definition we have, for all $n \in \{ 1, \dots, N^d\}$, 
    \begin{align*}
        \operatorname{dist}(\vec{e}(n), L)^2 = \sum_{l = S}^M \left| \langle \vec{e}(n), \nu_l \rangle \right|^2 = \sum_{l = S}^M \left| \sum_{j = 1}^M e_j(n) \ \overline{\nu_{j,l}} \right| ^2 
        = \sum_{l = S}^M  \sum_{j,k = 1}^M \overline{\nu_{j,l}} e_j(n)  \overline{e_k(n)} \nu_{k,l},  
    \end{align*}
    and thus 
    \begin{align*}
        \sum_{n = 1}^{N^d} \operatorname{dist} \left (\vec{e}(n), L \right)^2  = \sum_{n = 1}^{N^d} \sum_{l = S}^M  \sum_{j,k = 1}^M \overline{\nu_{j,l}} e_j(n)  \overline{e_k(n)} \nu_{k,l}  = \sum_{l=S}^M \left\langle \mathcal{E}_{N^d} \nu_l,\nu_l \right\rangle = \operatorname{tr}((I-\pi_L)\mathcal{E}_{N^d}).
    \end{align*}
    By estimating the summation, we obtain 
    \begin{equation} \label{4.8}
        \operatorname{tr}((I-\pi_L)\mathcal{E}_{N^d}) = \sum_{n = 1}^{N^d} \operatorname{dist} \left(\vec{e}(n), L\right)^2  \leq N^d\max _{n = 1, \dots, N^d} \operatorname{dist} \left(\vec{e}(n),L \right)^2.
    \end{equation}
    Thus, we can find an $n \in \{1, \dots, N^d \}$ which satisfies \eqref{usef diseq lem}.
    \end{proof}
    We let $0 \leq \varepsilon_1 \leq \varepsilon_2 \leq \ldots \leq \varepsilon_M$ denote the eigenvalues of $\mathcal{E}_{N^d}$. Applying the mini-max principle we obtain
    \begin{equation} \label{inf of the trace}
        \inf_{\dim L = S - 1} \operatorname{tr} \left( (1 - \pi_L) \mathcal{E}_{N^d} \right) = \varepsilon_1 + \varepsilon_2 + \ldots + \varepsilon_{M - S + 1} =: E_S.
    \end{equation} 
    \begin{oss}
        If we assume $
        \{e_1, \ldots, e_M\}$ to be an orthonormal family in $  \mathbb{C}^{N^d}$, then $\mathcal{E}_{N^d} = I$ and $\varepsilon_i =1$ for $i =1,\dots, M$. Thus, from \eqref{usef diseq lem} and \eqref{inf of the trace} we obtain
        \begin{equation*}
            \max_{n = 1, \dots N^d} \operatorname{dist}\left( \vec{e}(n), L \right)^2 \geq \frac{M - S + 1}{N^d}.
        \end{equation*}
    \end{oss}
    Now let $\{e_1, \dots, e_M\} $ and  $\{f_1, \dots, f_M\}$ be two linearly independent families in $\mathbb{C}^{N^d}$, and consider the notation \eqref{notation for vector e} for $\vec{e}(n)$ and $\vec{f}(n)$. Applying Lemma \ref{lem sjostr} together with \eqref{inf of the trace} to the first family, we successively choose indices 
    \begin{equation} \label{indices}
        n_1, \ldots, n_M \in \{ 1, \dots, N^d\}, 
    \end{equation}
    such that
    \begin{align*}
        &\left\|\vec{e}(n_1) \right\|^2   \geq N^{-d}\operatorname{tr} \,((I - \pi_{\{0\}})\mathcal{E}_{N^d}) \geq N^{-d} \inf_{\operatorname{dim} L = 0}\operatorname{tr} \,((I - \pi_L)\mathcal{E}_{N^d}) = N^{-d} E_1, \\
        &\operatorname{dist} \left( \vec{e}(n_S), \mathbb{C} \vec{e}(n_1) \oplus \cdots \oplus \mathbb{C} \vec{e}(n_{S-1}) \right)^2 \geq N^{-d} \inf_{\operatorname{dim} L = S-1}\operatorname{tr} \,((I - \pi_L)\mathcal{E}_{N^d}) = N^{-d} E_S,
    \end{align*}
    for $s = 2, \dots, M$.
    Let $ \{\nu_1, \nu_2, \ldots, \nu_M \}$ be the Gram-Schmidt orthonormalization of the basis
    $\{ \vec{e}(n_j) \}_{j=1}^M$,
    so that
    \begin{equation}\label{estimate con cj}
    \begin{aligned}
        &\vec{e}(n_1) = c_1 \nu_1 , \quad \text{where } |c_1| \geq \left( N^{-d} E_1 \right) ^{\frac{1}{2}}, \\
        &\vec{e}(n_S) \equiv c_S \nu_S \bmod (\nu_1, \ldots, \nu_{S-1}), \quad \text{where } |c_S| \geq \left( N^{-d} E_S\right) ^{\frac{1}{2}}, \quad \text{for } S = 2, \dots, M.
    \end{aligned}
    \end{equation}
    
    We consider the matrix \( E \) defined as in \eqref{def of E and F}. Expressing the vectors \( \vec{e}(n_1), \vec{e}(n_2), \ldots, \vec{e}(n_M) \) in the basis \( \nu_1, \ldots, \nu_M \) does not affect the absolute value of its determinant and hence $ |\det E| = |c_1 \cdots c_M|$. It follows from \eqref{estimate con cj} that
    \begin{equation}\label{formula for det}
        |\det E| \geq \frac{(E_1 E_2 \cdots E_M)^{1/2}}{(N^d)^{M/2}}.
    \end{equation}
    We now consider the matrix \( B \) defined as in \eqref{def of M} for the indices chosen in \eqref{indices}. Supposing that 
    \begin{equation} \label{hypoth on F and E} 
    f_j := \overline{e_j}, \ j=1, \dots, M,
    \end{equation}
    we have $F^* = E^t$, and consequently 
    \begin{equation} \label{definitio of B} 
        B = E E^t.
    \end{equation}
    Hence, under these hypotheses, from \eqref{formula for det}, we get
    \begin{equation} \label{5.16}
        |\det B| \geq \frac{E_1 E_2 \cdots E_M}{N^{Md}}.
    \end{equation}
    We want to estimate the singular values $s_1(B) \geq \ldots \geq s_M(B)$ of the matrix $B$. Writing $s_j := s_j(B)$, we have, for every $k = \{ 1, \dots, M \}$
    \begin{equation} \label{5.18}
        s_1^M  \geq s_1^{k-1} s_k^{M-k+1} \geq \prod_{j=1}^M s_j = | \det B |.
    \end{equation}
    We recall that $s_1 = \| B \|$, and combining \eqref{5.16} and \eqref{5.18}, we obtain
    \begin{equation*}
        s_1 \geq  \frac{(E_1 E_2 \cdots E_M) ^{\frac{1}{M}}}{N^{d}}.
    \end{equation*}
    and for each $k \in \{1,\dots,M\}$,
    \begin{align*}
            s_k \geq& \left( \left(\frac{1}{s_1}\right)^{k-1}  \prod_{j = 1}^M \left( \frac{E_j}{N^d}\right) \right) ^{\frac{1}{M-k+1}} \! \! \!= \left( s_1^{M- k+ 1}  \prod_{j = 1}^M \! \left( \frac{E_j}{s_1N^d}\right) \right) ^{\frac{1}{M-k+1}} \! \! \!= s_1 \left( \prod_{j = 1}^M \! \left( \frac{E_j}{s_1 N^d}\right) \right)^{\frac{1}{M-k+1}}.
    \end{align*}
    Summing up, we have proved the following.
    \begin{prop} \label{4.1.3}
        Under the assumption \eqref{hypoth on F and E}, the following two inequalities are true
        \begin{align}
            &s_1 \geq  \frac{(E_1 E_2 \cdots E_M) ^{\frac{1}{M}}}{N^{d}}, \label{est s1}\\
            &s_k \geq s_1 \left(   \prod_{j = 1}^M \left( \frac{E_j}{s_1N^d}\right) \right)^{\frac{1}{M-k+1}}, \quad \text{for } k = 1, \dots, M. \label{est s2}
        \end{align}
    \end{prop}
    We express the matrix $B$ as a matrix product involving a diagonal potential. Considering the indices $n_1, \dots, n_M \in \{ 1, \dots, N^d \}$ chosen above and recalling \eqref{entries of M}, we obtain for every $l,m = 1, \dots, M$,
    \begin{equation} \label{decomp of B}
    (B)_{l,m} =  \sum_{j=1}^M  e_l(n_j) \overline{f_m(n_j)} = \sum_{n=1}^{N^d}  e_l(n) \overline{f_m(n)} q(n),
    \end{equation}
    where $q \in \{ 0,1 \}^{N^d}$ is defined as
    \begin{equation} \label{def of q}
        q(n) := 
        \begin{cases} 
            1, \quad \text{if } \ \exists \ j : n = n_j\\
            0, \quad \text{otherwise}
        \end{cases},
    \end{equation}
    for every $n = 1, \dots, N^d$. In particular, $q$ has $M$ non-zero entries. 
    
    \subsection{Improving singular value bounds via potential perturbations} \label{section 3.2}
    In this section we adapt the work of Sj\"ostrand in \cite[Chapters 5–6]{Sjostrand 2} to the discrete case. 
    From now on, we adopt the same setting and notation as in Section \ref{Chapter 2}.  
    We set
    $  h = (2 \pi N)^{-1}, \  N \in \mathbb{N}^*$  and $ N^{-1} \ll \alpha \ll 1$ and let
    \begin{equation} \label{asymptotic expansion of the symbol}
        p \in C^\infty(\mathbb{T}^{2d}) \quad \text{admit the asymptotic expansion} \quad  p \sim p_0 + h p_1 + \dots \quad \text{in } S(1,1),
    \end{equation}
    We consider $P := p_N$ the restriction of its Weyl quantization to $\mathcal{H}_h^d$ as in \eqref{def of P_N,alpha}, so that $ P : \mathbb{C}^{N^d} \rightarrow \mathbb{C}^{N^d}.$ 
    We fix \( z \in \mathbb{C} \) and assume the symmetry condition 
    \begin{equation} \label{symmetry assumtpion 0}
        p(x,\xi) = p(x,-\xi), 
    \end{equation}
    which implies that
    \begin{equation} \label{symmetry assumption}
        P^* = \Gamma P \Gamma,
    \end{equation}
    where \( \Gamma u := \overline{u} \) denotes the antilinear operator of complex conjugation.  
    We observe that condition \eqref{symmetry assumption} remains valid under the addition of any diagonal matrix to \( P \). We assume that there exists \( \kappa \in ]0,1] \) such that
    \begin{equation} \label{bounded volume of the operator}
        V_z(t) := \lambda \left( \left\{ \rho \in \mathbb{T}^{2d} \; ; \; |p_0(\rho) - z|^2 \leq t \right\} \right) = \mathcal{O}(t^{\kappa}), \quad 0 \leq t \ll 1.
    \end{equation}

    \medskip
    Under these assumptions, Proposition \ref{teorema 2.2.2} yields the following.
    \begin{cor} \label{num of eigenval}
         For $ N^{-1} \ll \alpha \ll 1$, if \eqref{bounded volume of the operator} holds, the number $M(\alpha)$ of eigenvalues of $(P - z)^*(P - z)$ in $[0, \alpha]$ satisfies $$ M(\alpha) = \mathcal{O}(\alpha^{\kappa} N^{d}).$$
    \end{cor} 
    Let $q_0 \in \mathbb{C}^{ N^d}$ and $V_{q_0} := \operatorname{diag}(q_0(j); \ j= 1, \dots, N^d)$ its associated diagonal matrix. Let $\delta_0 \geq 0$ be a small parameter such that 
    \begin{equation} \label{hyp per perturb suff piccola}
    \| \delta_0V_{q_0} \| = \delta_0 \|  q_0 \|_\infty\ll N^{-1}, 
    \end{equation} 
    and consider the perturbed operator $P_0 := P + \delta_0 V_{q_0}$. We observe that under the hypothesis \eqref{hyp per perturb suff piccola}, by applying the mini-max theorem, Corollary \ref{num of eigenval} still applies after replacing $P$ with $P_0$.
    \\
    
    We fix $\tau_0 \in ]0,(CN^{-1})^{1/2}]$ for a new constant $C>0$, and denote by  
    $ 0 \leq t_1(P_0 - z) \leq \ldots \leq t_M(P_0 - z) $, all the singular values of $P_0 - z$ lying in the interval $[0,\tau_0[$.
    Following the previous notation, we have $M = M_0(\tau_0^2)$. Thus, if we take $  \alpha = CN^{-1}$, we obtain 
    \begin{equation} \label{estimate on M}
    M \leq  M_0(\alpha) = \mathcal{O}(\alpha^{ \kappa} N^{d}) = \mathcal{O} (N^{d - \kappa}).
    \end{equation}
    We recall that \( M \) is defined by the condition
    \begin{equation} \label{hyp for t_M}
        t_M(P_0 - z) < \tau_0 \leq t_{M+1}(P_0 - z).
    \end{equation}
    We fix \( \theta \in \left] 0, \frac{1}{4} \right[ \) and \( \epsilon_0 > 0 \) and obtain the following.
    \begin{prop} \label{4.2.2}
        Under the assumptions \eqref{asymptotic expansion of the symbol}, \eqref{symmetry assumtpion 0}, \eqref{bounded volume of the operator} and \eqref{hyp per perturb suff piccola} and supposing $N$ is large enough, the following statements hold:
        \begin{enumerate}
            \item Suppose \( M \) is sufficiently large. Then there exists a vector \( q \in \{0,1\}^{N^d} \), and the corresponding potential \( V_q := \operatorname{diag}(q(n); \ n = 1, \dots, N^d) \) with $\| V_q \| \leq 1$ such that, setting
            $$
                P^{\delta} := P_0 + \delta V_q, \qquad \delta := \frac{\tau_0}{C} N^{-2d},
            $$
            for some sufficiently large constant \( C > 0 \), the following estimates hold:
            \begin{align}
                t_{\nu}(P^{\delta} - z) &\geq t_{\nu}(P_0 - z) - \frac{\tau_0 N^{-d}}{C} 
                \geq \left( 1 - \frac{N^{-d}}{C} \right) t_{\nu}(P_0 - z), \quad \text{for } \nu > M, \label{first eq} \\
                t_{\nu}(P^{\delta} - z) &\geq \tau_0 N^{-M_1}, \quad \text{for } \ [M(1 - \theta)] + 1 \leq \nu \leq M, \label{second eq}
            \end{align}
            where
            \begin{equation} \label{N2}
                M_1 := 4d + \epsilon_0,
            \end{equation}
            and \( [a] := \max(\mathbb{Z} \cap (-\infty, a]) \) denotes the integer part of the real number \( a \).
    
            \item If \( M = \mathcal{O}(1) \), the same conclusion holds, provided that estimate~\eqref{second eq} is replaced by
            \begin{equation} \label{Caso M piccolo}
                t_M(P^{\delta} - z) \geq \tau_0 N^{-M_1}.
            \end{equation}
        \end{enumerate}
    \end{prop}

    \begin{proof}
        0. First notice that if $M=0$, then $t_1(P_0 -z) \geq \tau_0$ and the statement of the Proposition follows with $q = 0$. Hence, carrying on, we may assume that $M \geq$ 1.
        \\
        
        \noindent 1. We suppose $M$ sufficiently large, i.e. $M \gg 1$. Let $  \{e_1, \dots, e_M\} $ be an orthonormal family of eigenvectors of the operator $(P_0 - z)^*(P_0 - z)$, associated with its first $M$ eigenvalues $\{ t_1^2, \dots t_M^2\}$, analogously to those we introduced in Section \ref{Grushin problem for the unperturbed operator} for the unperturbed operator \( P \). By the symmetry assumption \eqref{symmetry assumption}, a corresponding orthonormal family of eigenvectors of $(P_0 - z)(P_0-z)^*$ associated to the same eigenvalues, is given by 
        \begin{equation} \label{def of f tilde}
        \widetilde{f_j} := \Gamma e_j = \overline{e_j}, \quad j= 1, \dots, M.
        \end{equation}
        We know that the  orthonormal family $\{f_1, \dots, f_M \}$ defined as in \eqref{def of f_i} (with $P$ replaced by $P_0$), and $\{ \widetilde{f_1}, \dots, \widetilde{f_M} \}$ generate the same vector space. 
        We consider a Grushin problem for \( P_0 - z \), following similar steps as in Section \ref{Grushin problem for the unperturbed operator}. In analogy with that construction, we define $R_+^0: \  \mathcal{H}_h^d \to \mathbb{C}^M, R_-^0: \mathbb{C}^M \to \mathcal{H}_h^d$  by
        \begin{align*} 
        R_+^0(u) = \sum_{i=1}^M \langle u, e_i \rangle \, \delta_i, \qquad R_-^0(u_-) = \sum_{i=1}^M u_-(i) \widetilde{f_i}.
        \end{align*}
        Then
        \begin{equation*} 
        \mathcal{P}^0 = \begin{pmatrix}
            P_0 - z & R_- \\
            R_+ & 0
        \end{pmatrix} :
        \mathcal{H}_h^d \oplus \mathbb{C}^M \longrightarrow \mathcal{H}_h^d \oplus \mathbb{C}^M,
        \end{equation*}
        is invertible, with bounded inverse
        \begin{equation}
        \mathcal{E}^0 =
        \begin{pmatrix}
            E^0 & E_+^0 \\
            E_-^0 & E_{-+}^0
        \end{pmatrix}.
        \end{equation}
        and the corresponding estimates in \eqref{diseq matr} remain valid.
    
        Since the properties stated in \eqref{properties of the eigenvectors e_j and f_j} do not necessarily hold in this setting, we cannot assert that \( E_{-+}^0 = \operatorname{diag}( t_j; \, j = 1,\dots,M) \). However, we know that $ t_j(E_{-+}^0) = t_j(P_0 - z), \  j = 1, \dots, M$.
        Since $z$ is fixed, we may assume without loss of generality that \( z = 0 \).
        \\
        
        \noindent \textsc{Case 1:} If $ s_j(E_{-+}^0) \geq \tau_0 N^{-M_1},  \text{for } 1 \leq j \leq M - [(1- \theta)M]$, the proposition is proved by taking the trivial perturbation $q= \boldsymbol{0} \in \mathbb{C}^{N^d}$.
        \\
        
        \noindent \textsc{Case 2:} We suppose there exists $j, \ 1 \leq j \leq M - [(1- \theta)M]$ such that
        \begin{equation} \label{Case 2}
        s_j(E_{-+}^0) < \tau_0 N^{-M_1}. 
        \end{equation}
        Since the family $\{\widetilde{f_j}\}$ satisfies hypothesis \eqref{hypoth on F and E}, we may apply the results of Section \ref{Section from Sjostrand} within the present setting and notation. It follows from \eqref{def of E and F}, \eqref{definitio of B}, \eqref{decomp of B} and \eqref{def of q} that there exists a matrix \( E \in \mathbb{C}^{M\times M} \), depending on the eigenvectors $\{ e_i\}_{i=1,\dots,M}$, and a vector \( q \in \{0,1\}^{N^d} \), with exactly $M$ non-zero entries equal to $1$, such that 
        \begin{equation} \label{forma esplicita potenziale}
            E E^t = E_-^0 V_q E_+^0 
        \end{equation}
        where
        \begin{equation} \label{potential}
            V_q := \operatorname{diag}\left( q(n);\ n=1, \dots N^d\right) \in \mathbb{C}^{N^d \times N^d}
        \end{equation}
        has the form of a potential. We note that $V_q$ is entirely determined by the choice of the indices $n_1, \dots, n_M$ in \eqref{indices}, and that
         \begin{equation} \label{norm of Vq 0}
            \| V_q \| = \| q \|_{\infty} = 1.
        \end{equation} 
        The matrix \( V_q \) is the perturbation we apply to the operator \( P_0 \) to derive the conclusions of the theorem. Accordingly, for $ \delta \geq 0$, we define the perturbed operator $ P^\delta := P_0 + \delta V_q$.
       
        We now set up a new Grushin problem for the perturbed operator, following the same procedure as in Section \ref{Grushin for perturbed op}.  
        We therefore consider the matrix
        \begin{equation} \label{perturbed Grushin}
        \mathcal{P}^\delta :=\begin{pmatrix}
                P^\delta & R_-^0 \\
                R_+^0    & 0
            \end{pmatrix},
        \end{equation}
        and we observe that condition \eqref{diseq Q} holds in our framework if and only if $ \delta \leq  \frac{\tau_0}{2}$.
        Hence, taking 
        \begin{equation} \label{hyp on delta}
            0 \leq \delta \leq \frac{\tau_0}{2},
        \end{equation} the Grushin problem \eqref{perturbed Grushin} is well-posed, and its inverse has the form
        \begin{equation} \label{def of E cal}
            \mathcal{E}^\delta =
            \begin{pmatrix}
                E^\delta      & E_+^\delta \\
                E_-^\delta    & E_{-+}^\delta
            \end{pmatrix},
        \end{equation}
        where the blocks are defined analogously to \eqref{Sviluppo di Neumann di cal(E)} and they satisfy the same norm estimates.
        
        Our next goal is to estimate the singular values of $P^\delta$. Recalling Lemma \ref{estimates for sing values} and applying it to the perturbed Grushin system \eqref{perturbed Grushin}, we obtain by \eqref{estimate for E^delta}, \eqref{estimate for E+-^delta}, \eqref{estimate for E-+^delta}, \eqref{diseq matr},
        \begin{equation} \label{result sing val of Pdelta and E-+delta}
            \frac{1}{\widetilde{C}} t_\nu (E_{-+}^\delta) \leq t_\nu (P^\delta)\leq \widetilde{C} t_\nu (E_{-+}^\delta)  , \qquad  \nu = 1, 2, \dots, M,
        \end{equation}
        for a constant $\widetilde{C}>0$.
        From \eqref{Sviluppo di Neumann di cal(E)}, the lower right block of \eqref{def of E cal} is defined as
        \begin{equation} \label{4.42}
            E_{-+}^\delta 
            = E_{-+}^0 - \delta E_-^0 V_q E_+^0  
              + \sum_{k=2}^{+\infty} (-1)^k 
                E_-^0 (\delta V_q E^0)^{k-1} (\delta V_q) E_+^0.
        \end{equation}
        Applying the estimates \eqref{diseq matr} with the hypothesis \eqref{hyp for t_M} to our setting yields
        \begin{equation} \label{new estimate for E}
            \| E_\pm^0 \| = 1, \qquad \| E^0 \| \leq \frac{1}{t_{M+1}(P_0)} \leq \frac{1}{\tau_0}.
        \end{equation}
        We fix $ 1 \leq \nu \leq M - [ (1-\theta) M ]$ and we take the index $j$ defined in the assumption \eqref{Case 2}, so that $\nu + j - 1 \leq M$. Applying Ky Fan's inequalities (\cite[Corollary 2.2]{Gohberg}) to the decomposition \eqref{4.42} yields
        \begin{align} \label{sum of sing values}
            s_{\nu}(E^\delta_{-+}) \geq  s_{\nu+j-1}(-\delta E^{0}_{-} V_q E^{0}_{+}) - s_{j}(E^{0}_{-+}) 
            - s_1\left(\sum_{k=2}^{+\infty} (-1)^k E_-^0(\delta V_q E^0)^{k-1} (\delta V_q) E_+^0 \right). 
        \end{align}
        We now estimate each term separately. The second term can be estimated using hypothesis \eqref{Case 2}, which implies $ -s_j(E_{-+}^0) > -\tau_0 N^{-M_1}$. 
        For the third term, we recall \eqref{norm of Vq 0}, \eqref{new estimate for E}, obtaining
           \begin{align}
            s_1\left(\sum_{k=2}^{+\infty} (-1)^k E_-^0(\delta V_q E^0)^{k-1} (\delta V_q) E_+^0 \right)  
            \leq \sum_{k=2}^{+\infty}  \delta ^k  \| V_q \| ^k \| E^0 \|^{k-1} 
            \leq 2 \delta ^2 \tau_0^{-1}, \label{third term}  
        \end{align} 
        where we used the fact that
        $\delta \| V_q \| \| E^0 \| \leq \frac12$,
        ensuring the convergence of the series. 
        For the first term, from \eqref{def of f tilde} we can invoke Proposition \ref{4.1.3} and we observe that \eqref{new estimate for E} implies $s_1(E_-^0V_q E_+^0) \leq 1$. Thus, for $k=1,\dots,M$ we rewrite \eqref{est s2} as
        \begin{align} \label{firste useful diseq}
            s_{k}(E^0_-V_q E^0_+) \geq  s_1(E_-^0V_q E_+^0)^{-\frac{ k -1}{M - k +1}} \left(  {M!} \right)^{\frac{1}{M - k + 1}} N^{-\frac{Md}{M- k + 1}}  \geq \left(  {M!} \right)^{\frac{1}{M - k + 1}} N^{-\frac{Md}{M- k + 1}}. 
        \end{align}
        Now, if $1 \leq k \leq \frac{M}{2}$, then $\frac{M}{M-k+1} \leq 2$. Thus, taking $k=\nu + j - 1 \leq \frac{M}{2} $ for $\theta \in ]0,\frac{1}{4}[$, we get
        \begin{equation} \label{estimate for EF*}
            s_{\nu + j -1}( - \delta E_-^0 V_q E_+^0) \geq \delta \left(  {M!} \right)^{\frac{1}{M - k + 1}} N^{-\frac{Md}{M- k + 1}}  \geq \delta N^{-2d}.
        \end{equation}
        Combining \eqref{sum of sing values}, \eqref{Case 2}, \eqref{third term}, and \eqref{estimate for EF*} we obtain the following inequality
        \begin{align} \label{almost final estimate}
            s_{\nu}(E^\delta_{-+}) \geq \delta N^{-2d} - \tau_0 N ^{-M_1} - 2 \delta^2 \tau_0^{-1}  , \quad 
            1\leq \nu\leq M - \left[(1-\theta)M \right].
        \end{align}
        We fix $1\leq \nu\leq M - \left[(1-\theta)M \right]$ and we choose
        \begin{equation} \label{def of delta}
            \delta = \frac{\tau_0}{C}N^{-2d} ,
        \end{equation}
        where $C>0$ is a sufficiently large constant. Taking $N$ sufficiently large, condition \eqref{hyp on delta} on $\delta $ is respected and substituting this value in \eqref{almost final estimate}, we obtain
        \begin{equation} \label{final estimate}
             s_{\nu}(E^\delta_{-+}) \geq \delta N^{-2d} \left( \frac{1}{e^{1/2}} - C N ^{-\epsilon_0} - \frac{2}{C} \right) \geq \frac{\delta}{C_0}N^{-2d},
             \end{equation}
        for a new constant $C_0>0$. Thus, by possibly taking a larger $N$, we obtain $s_{\nu}(E^\delta_{-+})\geq \widetilde{C} \tau_{0} N^{-M_1}$, with the same constant $\widetilde{C}$ that appears in \eqref{result sing val of Pdelta and E-+delta} and $M_1$ defined in \eqref{N2}. It follows that $$t_\nu (P^\delta)  \geq\tau_0 N^{-M_1}, \quad \text{for }  1+ \left[(1-\theta) M \right]\leq \nu \leq M,$$
        and therefore \eqref{second eq} is proved.
        On the other hand, if $\nu > M$, Ky Fan's inequalities and \eqref{norm of Vq 0} yield
        \begin{align*}
            t_\nu (P^\delta)\geq t_\nu (P_0) - s_1(\delta V_q) 
            = t_\nu (P_0) -\frac{\tau_0 N^{-2d}}{C} \geq t_\nu(P_0) \left( 1 - \frac{ N^{-d}}{C}\right).
        \end{align*} 

        \noindent
        2. Finally, we consider the case where $M = \mathcal{O}(1)$. The argument proceeds as before, except that it now reduces to two subcases concerning only the $M$-th singular value of $E_{-+}^0$. 
        Carrying out the proof in this way yields inequality \eqref{final estimate} with $\nu = 1$, thereby establishing \eqref{Caso M piccolo}.        
    \end{proof}
    \begin{oss} \label{remark for the iteration}
        From the preceding proof, we note that Proposition~\ref{4.2.2} remains valid provided that $M$ satisfies only the condition $t_{M+1}(P_0 - z) \;\geq\; \tau_0$.
    \end{oss}

    \subsection{Iterative construction of a potential for controlling the smallest singular value}
    We now adopt the same hypotheses considered at the beginning of the previous section and we iterate the construction of the potential considered in the proof of Proposition \ref{4.2.2}, starting from the unperturbed operator $P$. Our goal is to determine a new potential $V_{q_0}$, with $q_0 \in \mathbb{C}^{N^d}$, and establish a lower bound for the smallest singular value of the perturbed operator $P + \delta V_{q_0} - z,  z \in \mathbb{C}$, where $ \delta \geq 0$. 
    Let $p$ be as in \eqref{asymptotic expansion of the symbol} and \eqref{symmetry assumtpion 0}. We fix $z \in \mathbb{C}$ that satisfies \eqref{bounded volume of the operator} for some $\kappa \in ]0,1]$, and set $P^{(0)} := P$.
    Let $\tau_0^{(0)} \in ]0, (CN^{-1})^{1/2}]$, for some constant $C > 0$, and denote by 
    \begin{equation} \label{def of M(0)}
    M^{(0)}= M^{(0)} (\tau_0^{(0)}),
    \end{equation}
    the number of singular values of $P^{(0)} - z$ lying in the interval $[0, \tau_0^{(0)}[$. If \(M^{(0)}=0\), we take \(q_0=0\), since \(t_1(P-z)\geq\tau_0^{(0)}\); no perturbation is needed. Henceforth, we assume that \(M^{(0)}\geq1\). First, we assume that $M^{(0)}$ is sufficiently large and we apply Proposition \ref{4.2.2} in this framework. Consequently, there exists a vector $q^{(1)} \in \{0,1\}^{N^d}$ with at most $M^{(0)}$ non-zero entries, defining the potential $V_{q^{(1)}}$ with $\| V_{q^{(1)}} \| \leq 1$, and a coefficient 
    $\delta^{(1)} := C^{-1}\tau_0^{(0)} N^{-2d}$, with a sufficiently large constant $C > 0$, such that the singular values of
    $ P^{(1)} := P^{(0)} + \delta^{(1)} V_{q^{(1)}}$
    satisfy \eqref{first eq} and \eqref{second eq} with the corresponding updated notations. In particular, setting
    $$
    \tau_0^{(1)} := \tau_0^{(0)} N^{-M_1}, 
    \quad M^{(1)} := [M^{(0)}(1-\theta)],
    $$
    with $M_1$ defined as in \eqref{N2} and $\theta \in ]0,\frac{1}{4}[$, we obtain
    $
    t_{\nu}(P^{(1)} - z) \geq \tau_0^{(1)}, 
    \ \text{for }  M^{(1)} +1 \leq \nu \leq M^{(0)}
    $.
    Thanks to Remark \ref{remark for the iteration}, we can apply Proposition \ref{4.2.2} once more, replacing $ (P^{(0)}, M^{(0)}, \tau_0^{(0)})$ with $ (P^{(1)}, M^{(1)}, \tau_0^{(1)}),$ while keeping the same value of $M_1$. Since
    $$
    \| \delta^{(1)} V_{q^{(1)}} \| \leq \frac{\tau_0^{(0)}}{C} N^{-2d} \ll N^{-1},
    $$
    condition \eqref{hyp per perturb suff piccola} holds. Thus, we iterate this procedure, applying at each $k$-th step the first case of Proposition \ref{4.2.2}, with the initial conditions given by the outcome of the previous step. In this way, we construct a sequence $(P^{(k)}, M^{(k)}, \tau_0^{(k)}), \ k = 0, 1, \dots, k(M)$,
    where the final index $k(M)$ is chosen so that $M^{(k(M))}$ remains of the order of a large positive constant, ensuring that the first case of the proposition can still be applied. Therefore at each $(k+1)$-th step of the iteration, we obtain $q^{(k+1)} \in \{ 0,1\}^{N^d}$, with at most $M^{(k)}$ non-zero entries, such that
    \begin{align}
        &t_\nu (P^{(k+1)}-z) \geq \tau_0^{(k+1)}, \quad M^{(k+1)} + 1 \leq \nu \leq M^{(k)}, \label{4.59}\\
        & t_\nu (P^{(k+1)}-z) \geq t_\nu (P^{(k)}-z) - \frac{\tau_0^{(k)} N^{-d}}{C}, \quad \nu > M^{(k)}. \label{4.60}, 
    \end{align}
    with
    \begin{align} \label{def of M(k)}
        \tau_0^{(k+1)} := \tau_0^{(k)} N^{-M_1}, \qquad M^{(k+1)} := [M^{(k)}(1-\theta)], \qquad P^{(k+1)} := P^{(k)} + \delta^{(k+1)} V_{q^{(k+1)}},
    \end{align}
    and
    \begin{align}
         \delta^{(k+1)} = \frac{\tau_0^{(k)}}{C} N^{-2d},\quad \text{and } \ \| V_{q^{(k+1)}}\| \leq 1,  \label{delta(k)}
    \end{align}

    We observe that \eqref{4.59} guarantees that, at each iteration step, the assumption required in Proposition \ref{4.2.2} holds in the weaker form described in Remark \ref{remark for the iteration}. From the iteration, we also deduce that 
    \begin{equation} \label{tau_0 iter}
        \tau_0^{(k)} = \tau_0^{(0)}N^{-M_1 k},
    \end{equation}
    and therefore
    \begin{equation} \label{delta iterated}
        \delta^{(k+1)} = \frac{\tau_0^{(0)} N^{-M_1 k}}{C} N^{-2d} = \delta^{(1)} N^{-M_1 k}.
    \end{equation}
    Hence, it follows that
    \begin{equation*}
        P^{(k+1)} = P^{(k)} + \delta^{(k+1)} V_{q^{(k+1)}} = P^{(0)} + \sum_{j = 1}^{k+1} \delta^{(j)} V_{q^{(j)}} =  P + \delta^{(1)} \sum_{j = 0}^{k} N^{-M_1 j} V_{q^{(j+1)}}.
    \end{equation*}
    Taking the norm of the perturbative term, we find by \eqref{delta(k)}
    \begin{align*}
        \left\|  \delta^{(1)} \sum_{j = 0}^{k} N^{-M_1 j} V_{q^{(j+1)}}\right\| \leq \delta^{(1)} \sum_{j = 0}^{k} N^{-M_1 j}\|    V_{q^{(j+1)}} \| < 2 \delta^{(1)} \ll N^{-1}.
    \end{align*}
    This bound shows that the hypothesis \eqref{hyp per perturb suff piccola} is satisfied for every $k$.  
    Therefore, the application of Proposition \ref{4.2.2} at each step of the iteration is justified. 
    We observe that $M^{(k)}$ decays exponentially fast in $k$. In particular,
    \begin{equation} \label{exp decay of M}
        M^{(k)} \leq (1 - \theta)^k M^{(0)}.
    \end{equation}
    We iterate until $M^{(k(M))}$ is of the order of a large constant. In particular we require
    $$(1-\theta)^{k(M)} M^{(0)} \geq C,$$
    for some constant $ C \gg 1$, which is equivalent to
    \begin{equation} \label{def of k_0}
        k_0 := k(M) \leq \frac{\ln \left( M^{(0)}/C\right)}{\ln \left( 1/(1- \theta) \right)}.
    \end{equation}
    Beyond $k=k_0$ (or immediately if $M^{(0)}$ is bounded), we continue the iteration finitely many times by decreasing $M^{(k)}$ by one at each step until $1$ is reached. For these steps, we apply \eqref{Caso M piccolo} and obtain the same estimates. Given $\nu > M^{(0)} > \ldots > M^{(k)}$, we use \eqref{tau_0 iter} and iterate \eqref{4.60} to obtain
    \begin{align*}
        t_\nu (P^{(k)}-z) \geq t_\nu(P^{(0)}-z) - \tau_0^{(0)} \frac{N^{-d}}{C}\left(\sum_{j=0}^{k-1} N^{-jM_1} \right) \geq  t_\nu(P^{(0)} - z) - \tau_0^{(0)} \mathcal{O} \left(\frac{N^{-d}}{C} \right).
    \end{align*}
    Moreover, for $1 \leq \nu \leq M^{(0)}$, let $l = l(M^{(0)}) $ denote the unique iteration step such that 
    $ M^{(l)} + 1 \leq \nu \leq M^{(l-1)}$,
    which is $\nu = M^{(l-1)}$, whenever $\nu \leq \widetilde{C}$.
    Then, applying \eqref{4.59} or \eqref{Caso M piccolo}, we obtain
    $$
        t_\nu (P^{(l)} -z ) \geq \tau_0^{(l)}.
    $$
    Hence, if $k >l$, from \eqref{4.59} and \eqref{4.60} we get
    \begin{align}
        t_\nu (P^{(k)}-z) \geq  t_\nu(P^{(l)}-z) -  \frac{N^{-d} \tau_0^{(l)}}{C} \left( \sum_{j=0}^{k-l-1} N^{-jM_1} \right) \geq \tau_0^{(l)} \left( 1- \mathcal{O} \left(\frac{N^{-d}}{C} \right)\right). \label{4.76}
    \end{align}
    Furthermore, since \(M^{(k_1)}=1\), estimate \eqref{Caso M piccolo} gives
    $$t_1(P^\delta -z) = t_{M^{(k_1)}} (P^{{(k_1+1)}} - z ).$$
    We now collect the perturbations constructed at the successive steps. Using
    \eqref{delta iterated}, we obtain
    \begin{equation*}
        \delta V_{q_0} := \sum_{k=1}^{k_1+1} \delta^{(k)} V_{q^{(k)}} =  \sum_{k=1}^{k_1+1} \delta^{(1)} N^{-M_1 (k-1)} V_{q^{(k)}} = \delta^{(1)}  \sum_{k=0}^{k_1} N^{-M_1 k} V_{q^{(k+1)}},
    \end{equation*}
    Accordingly, we define
        $$ q_0 :=  \sum_{k=0}^{k_1} N^{-M_1 k} q^{(k+1)} \in \mathbb{C}^{N^d}.$$
    Since \(q^{(k+1)}\) has at most \(M^{(k)}\) nonzero entries, each equal to \(1\), we have
    \begin{equation} \label{eq non capita}
    \begin{aligned}
        \|q_0\|_2 \leq \sum_{k=0}^{k_1} N^{-M_1k}\sqrt{M^{(k)}} \leq \sqrt{M^{(0)}}\sum_{k=0}^{k_1} N^{-M_1k} \leq \frac{\sqrt{M^{(0)}}}{1-N^{-M_1}} = \mathcal{O}\left(N^{\frac{d-\kappa}{2}}\right).
    \end{aligned}
    \end{equation}
    and similarly $ \| q_0 \|_{\infty}  < 2$ for $N$ sufficiently large.
    
    We have thus constructed a single perturbation \(V_{q_0}\) that encodes all
    the steps of the iteration. The preceding estimates yield the following
    proposition.
    \begin{prop} \label{4.2.3}
        Let $P$ defined as in Section \ref{section 3.2} and fix $z \in \mathbb{C}$ that satisfies \eqref{bounded volume of the operator}. Let $M_1 = 4 d + \epsilon_0,$  where \( \epsilon_0 > 0 \).
        Fix \( 0 < \tau_0 \leq \sqrt{C_\tau N^{-1}} \), $C_\tau>0$ and let \( M^{(0)} = \mathcal{O}(N^{d - \kappa }) \) be the number of singular values of \( P - z \) in \( [0, \tau_0[ \). Let \( 0 < \theta < \frac{1}{4} \) and let \( M(\theta) \gg 1 \) be a sufficiently large constant.
        
        We set \( M^{(k)},\ 1 \leq k \leq k_1 \) iteratively as follows. As long as \( M^{(k)} \geq M(\theta) \), we set \( M^{(k+1)} = [(1 - \theta)M^{(k)}] \). Let \( k_0 \geq 0 \) be the last \( k \) value obtained in this way. For \( k > k_0 \) we put \( M^{(k+1)} = M^{(k)} - 1 \), until we reach the value \( k_1 \) such that \( M^{(k_1)} = 1 \).
        We set $  \tau_0^{(k)} = \tau_0 N^{-k M_1}$, for $ \ 1 \leq k \leq k_1 + 1 $. 
        Then there exists $q_0 \in \mathbb{C}^{N^d}$ with its corresponding potential matrix \( V_{q_0}  \in \mathbb{C}^{N^d \times N^d} \), such that
        \begin{equation} \label{diseq for q0}
         \| q_0\|_{\infty} = \mathcal{O}(1),\qquad\| q_0\|_{2} = \mathcal{O}(N^{\frac{d - \kappa}{2}}).
        \end{equation}
        Moreover, if we define $ \delta = \frac{\tau_0}{C} N^{-2d}$, for $C> 0$ sufficiently large and \( P^\delta = P + \delta V_{q_0} \), then we have the following estimates on the singular values of \( P^\delta - z \):
        \begin{itemize}
          \item If $  \nu > M^{(0)} $,  
          \begin{equation} \label{1 result}
          t_\nu(P^\delta - z) \geq \left( 1 - \frac{N^{-d}}{C} \right ) t_\nu (P - z).
          \end{equation}
          \item If \( M^{(k)} + 1 \leq \nu \leq M^{(k-1)}, \) with \(\ 1 \leq k \leq k_1 \), then 
          \begin{equation} \label{2 result}
          t_\nu(P^\delta - z) \geq (1 - \mathcal{O}(N^{-d})) \tau_0^{(k)}.
          \end{equation}
          \item Finally, for \( \nu = M^{(k_1)} = 1 \), we have 
          \begin{equation} \label{3 result}
          t_1(P^\delta - z) \geq  \tau_0^{(k_1 + 1)} \geq (1 - \mathcal{O}(N^{-d})) \tau_0^{(k_1 + 1)}. \end{equation}
        \end{itemize}
    \end{prop}

    \subsubsection{Log-determinant estimates for the perturbed operator}
    We proceed under the same assumptions of Proposition \ref{4.2.3} and we consider the perturbed operator $ P^\delta(q_0) := P + \delta V_{q_0}$ arising from that result. Fixed $z \in \mathbb{C}$, our goal is to establish a lower bound for $\det(P^\delta(q_0)-z)$.
    We choose $C_0>C_\tau$ sufficiently large, independently of $N$, and consider
    \begin{equation} \label{def of M 2}
    M := M(C_0 N^{-1})= \mathcal{O}(N^{d - \kappa }),
    \end{equation}
    the number of singular values of $P^\delta(q_0) - z$ in the interval $[0,(C_0N^{-1})^{1/2}[$. In particular we are working with $\alpha = C_0N^{-1}$, so that the hypothesis of Proposition \ref{teorema  2.2.2} is satisfied. 

    We set up the $M$-dimensional Grushin problem for $P^\delta(q_0)-z$ as in Section \ref{Grushin for perturbed op}.
    Let $R>0$ satisfy \eqref{interval of definition of R}.
    By \eqref{diseq for q0} and the lower bound in
    \eqref{interval of definition of R}, choosing $C_2$
    sufficiently large ensures that
    \begin{equation} \label{bound on 2 norm of q0}
        \|q_0\|_2 < R/2.
    \end{equation}
    Moreover, recalling that $\alpha=C_0N^{-1}$ and
    $\delta=\tau_0N^{-2d}/C$, the upper bound in
    \eqref{interval of definition of R} implies $R \leq \frac{C\sqrt{C_0}}{6}\, \tau_0^{-1}N^{2d-\frac12},$. Thus, for every $q\in B_{N^d}(0,3R)$,
    \begin{equation} \label{estim for R}
        \delta\alpha^{-1/2}\|V_q\|
        = \delta\alpha^{-1/2}\|q\|_\infty
        \leq 3\delta\alpha^{-1/2}R
        \leq \frac12.
    \end{equation}
    Thus condition \eqref{diseq Q} holds throughout
    $B_{N^d}(0,3R)$ and, in particular, at $q=q_0$. Under this hypothesis the perturbed Grushin problem is well defined. We keep the same notations as in Section \ref{Grushin for perturbed op}, now emphasizing the dependence on $q_0$. In particular, we denote the matrices of the problem as $\mathcal{P}^\delta(q_0), R_\pm(q_0),$
    and the inverse matrix of $\mathcal{P}^\delta(q_0)$ is
    $$
    \mathcal{E}^\delta(q_0) = 
    \begin{pmatrix}
    E^\delta(q_0) & E_+^\delta(q_0) \\
    E_-^\delta(q_0) & E_{-+}^\delta(q_0)
    \end{pmatrix}.
    $$
    Our goal is to analyze the spectrum of the operator $P^\delta(q_0)$ by studying its characteristic polynomial. To this end, we recall formula \eqref{important sum for the final estimate}:
    \begin{equation} \label{formula per vari ln det}
        \ln |  \det(P^\delta(q_0)-z)| = \ln |  \det \mathcal{P}^\delta(q_0) |+ \ln |  \det E_{-+}^\delta(q_0)|, 
    \end{equation}
    and we study the two terms in the right-hand side separately.
    For the first term, by first applying \eqref{estimate per P cal} and then \eqref{formula per P cal}, we obtain
    \begin{align*}
        \ln|\det \mathcal{P}^\delta (q_0) | \geq   N^d \left( \int_{\mathbb{T}^{2d}} \ln \, |p_0(\rho) - z| d\rho - \mathcal{O}\left( N^{-\kappa} \ln(N) \right) \right)-  \mathcal{O}(\delta \alpha^{-1/2} \| V_{q_0} \|_{\operatorname{tr}}).
    \end{align*}
    Hence, by recalling that $\| V_{q_0} \|_{\operatorname{tr}} \leq N^{d/2} \| q_0 \|_{2}$ and $\|q_0\|_2 \leq R/2$, it follows from \eqref{estim for R} that
    \begin{equation} \label{estimate for P cal 3}
        \ln|\det \mathcal{P}^\delta(q_0)|
        \geq N^d \left( \int_{\mathbb{T}^{2d}}
            \ln|p_0(\rho)-z|\,d\rho - \mathcal{O}\left(N^{-\kappa}\ln N\right)
            - \mathcal{O}\left(N^{-d/2}\right) \right).
    \end{equation}
    On the other hand, for the second term we observe that
    \begin{equation} \label{determinant}
        |\det (E_{-+}^\delta (q_0)) | = \prod _{\nu = 1}^{M^{(0)}} t_\nu (E_{-+}^\delta (q_0)) \prod _{\nu = M^{(0)}+1}^{M} t_\nu (E_{-+}^\delta (q_0)),
    \end{equation}
    where $M$ is defined in \eqref{def of M 2} and $M^{(0)}$ in \eqref{def of M(0)}. 
    From \eqref{def of M(k)}, we deduce that, for $1 \leq k \leq k_0$,
    $ M^{(k-1)} - M^{(k)} \leq 1 + (1-\theta)^{k-1} \theta M^{(0)}$.
    We recall inequalities \eqref{result sing val of Pdelta and E-+delta} and, applying Proposition \ref{4.2.3}, we obtain for $1 \leq k \leq k_0$,
    \begin{align}
        \prod _{\nu = 1 + M^{(k)}} ^{M^{(k-1)}} t_\nu (E_{-+}^\delta (q_0)) &\geq \prod _{\nu = 1 + M^{(k)}} ^{M^{(k-1)}} \frac{1}{\widetilde{C}} (1 - \mathcal{O}(N^{-d})) \tau_0  N^{-k M_1} \geq  \left(\frac{1}{\widetilde{C}'} \tau_0 N^{-k M_1} \right) ^{1 + (1-\theta)^{k-1} \theta M^{(0)}}. \label{usef resu 2}
    \end{align}
    for a new constant ${\widetilde{C}'}>0$. On the other hand, if $k_0 < k \leq k_1$, then $M^{(k-1)} - M^{(k)} = 1$. Applying the same arguments, we obtain
    \begin{equation} \label{usef resu 3}
    t_{M^{(k-1)}} (E_{-+}^\delta (q_0)) \geq \frac{1}{\widetilde{C}'} \tau_0 N^{-k M_1}.
    \end{equation}
    Hence, by putting together \eqref{usef resu 2} and \eqref{usef resu 3}, we get
    \begin{align*}
        \ln \left( \prod _{\nu = 1}^{M^{(0)}} t_\nu (E_{-+}^\delta (q_0)) \right) \geq \sum_{k = 1}^{k_0} \ln   \left( \left( \frac{1}{\widetilde{C}'} \tau_0 N^{-k M_1} \right) ^{1 + (1-\theta)^{k-1} \theta M^{(0)}} \right) + \sum _{k = k_0 +1}^{k_1+1}\ln \left (  \frac{1}{\widetilde{C}'} \tau_0 N^{-k M_1}\right).
    \end{align*}
    Since $(1 + (1-\theta)^{k-1} \theta M^{(0)}) \geq 1$ for $k_0+1 \leq k \leq k_1 +1$, we obtain
    \begin{equation}\label{preresult chap 3}
        \ln \left( \prod _{\nu = 1}^{M^{(0)}} t_\nu (E_{-+}^\delta (q_0)) \right) \geq  - \sum_{k = 1}^{k_1 + 1} ({1 + (1-\theta)^{k-1} \theta M^{(0)}}) \left( \ln (\widetilde{C}') + \ln \left(\frac{1}{\tau_0} \right) +k M_1\ln (N) \right).
    \end{equation}
    We recall from \eqref{estimate on M} that $M^{(0)} = \mathcal{O}(N^{d - \kappa})$, and we consider the constant $M(\theta) \gg 1$ specified in the assumptions of Proposition \ref{4.2.3}. Hence, by \eqref{def of k_0}, we obtain
    $ k_1 + 1= \mathcal{O}( \ln (N))$.
    By applying this result to \eqref{preresult chap 3}, we obtain
    \begin{align}
        \ln \left( \prod _{\nu = 1}^{M^{(0)}} t_\nu (E_{-+}^\delta (q_0)) \right) 
        &\geq -  \left( \mathcal{O}(\ln (N)) +  \mathcal{O}(N^{d - \kappa })(1- (1-\theta)^{k_1+1}) \right) \notag\\
        &\hspace{0.5cm} \cdot \left( \ln (\widetilde{C}') + \ln \left(\frac{1}{\tau_0} \right) + \mathcal{O}(1) \left (\ln (N) \right ) ^2 \right). \label{sec fin res}
    \end{align}
    From \eqref{1 result}, the remaining factor in \eqref{determinant} can be bounded in the following way
    \begin{equation} \label{second bound sing values}
        \prod _{\nu = M^{(0)}+1}^{M} t_\nu (E_{-+}^\delta (q_0)) \geq \prod _{\nu = M^{(0)}+1}^{M} \frac{t_\nu (P-z)}{\widetilde{C}'}\geq  \left(\frac{\tau_0}{\widetilde{C}'} \right)^{M - M^{(0)}}
    \end{equation}
    Hence, applying the logarithm, we have
    \begin{equation} \label{terz fin res}
        \ln \left( \prod _{\nu = M^{(0)}+1}^{M} t_\nu (E_{-+}^\delta (q_0)) \right) \geq  -M \left(\ln\left( \frac{1}{\tau_0} \right) +\ln(\widetilde{C}') \right) \geq -\mathcal{O}(N^{d - \kappa})\ln\left( \left(\frac{1}{\tau_0} \right) + \ln (\widetilde{C}')  \right).
    \end{equation}
    Now we combine together \eqref{sec fin res} and \eqref{terz fin res} in order to find an estimate for $\ln|\det E_{-+}^\delta (q_0) |$. Since $\ln (\widetilde{C}')$ is a constant, 
    it is absorbed into the big-$\mathcal{O}$. Hence, there exists a constant $C_1 > 0$ such that
    \begin{equation} \label{final estimate ln det}
        \ln|\det E_{-+}^\delta (q_0) | \geq  - C_1 N^d \left( N^{-d} \ln(N) +  N^{-\kappa} \right) \left( \ln \left(\frac{1}{\tau_0} \right) +  \left (\ln(N) \right ) ^2 \right), 
    \end{equation}
    Finally, combining \eqref{formula per vari ln det}, \eqref{estimate for P cal 3}, and \eqref{final estimate ln det}, we obtain the following result.

    \begin{prop} \label{Prop for bound for q0}
    Under the assumptions of Proposition \ref{4.2.3}, for the potential $V_{q_0}$ constructed therein, we have that there exists $C>0$ such that
        \begin{align}
            &\ln|\det (P^\delta(q_0) - z) | - N^d\left(  \int_{\mathbb{T}^{2d}} \ln \, |p_0(\rho) - z| d\rho \right) \notag \\
            &\geq  - CN^d \left( N^{-\frac{d}{2}} + \left( N^{-d}\ln (N) +  N^{-\kappa} \right) \left( \ln \left(\frac{1}{\tau_0}\right) +  \left (\ln (N) \right ) ^2 \right) \right). \label{result end chapter 3}
        \end{align}
    \end{prop}
    
    \section{Perturbations by random potentials} \label{Chapter 4}
    In this section we first record some estimates derived from the computations in Section \ref{Chap 3}. Then we introduce a probability measure on \(\mathbb C^{N^d}\) which determine the perturbation potentials \(V_q\). Finally, we obtain a spectral estimate which holds in probability.
    
    Throughout this section we assume that \(p\in C^\infty(\mathbb T^{2d})\) satisfies \eqref{asymptotic expansion of the symbol} and \eqref{symmetry assumtpion 0}, so that its quantization $P := p_N$ satisfies \eqref{symmetry assumption}. We fix \(h=(2\pi N)^{-1}\) with \(N \in \mathbb{N}, \, N \gg 1 \) and $\tau_0 \in ]0, (C_\tau N^{-1})^{1/2}[$ where $C_\tau>0$. Moreover, we assume that \(z\in\mathbb C\) satisfies \eqref{bounded volume of the operator} for some \(\kappa\in ]0,1]\) and we recall \(q_0\in\mathbb C^{N^d}\) and \(\delta=\frac{\tau_0}{C}\,N^{-2d}\) (with \(C>0\) independent of \(N\)) from Proposition \ref{4.2.3}.
    We retain the notation introduced in Section \ref{Chap 3} for the associated Grushin problem matrices \(\mathcal P\) and \(\mathcal P^\delta\), with $\alpha = C_0N^{-1}$.

    \subsection{Deterministic log-determinant estimates for bounded perturbations} \label{section about estimate for a general potential} We consider $R>0$ satisfying \eqref{interval of definition of R}, so that 
    \begin{equation} \label{def of R}
    \| q_0 \|_{2} < \frac{R}{2}.
    \end{equation}
    In particular, we choose the constants in \eqref{interval of definition of R} so that condition \eqref{diseq Q} is still verified in the ball $B_{N^d}(0,3R)$.
    Let $q \in \mathbb{C}^{N^d}$ be an arbitrary vector such that $\| q \| _{2} < 3R.$ We observe that for its corresponding potential matrix $V_q$, the following estimate is true
    \begin{equation*}
        \| V_q \|_{\operatorname{tr}} \leq N^{d/2} \| q \|_{2} < \widetilde{C_1} N^{\frac{5d-1}{2}} \tau_0 ^{-1} ,
    \end{equation*}
    for a new constant $\widetilde{C_1} > 0$.
    Hence, by applying \eqref{estimate for mathcal(P)delta} with $\alpha = C_0 N^{-1}$, we obtain
     \begin{align}
        \ln|\det \mathcal{P}^\delta (q) |  \leq  N^d \left( \int_{\mathbb{T}^{2d}} \ln \, |p_0(\rho) - z| d\rho + \mathcal{O}\left( N^{-\kappa} \ln(N)\right) +   \mathcal{O} ( N^{-\frac{d}{2}} )\right),
    \end{align}
    where we used the fact that $\delta= \frac{\tau_0}{C}N^{-2d}$. Thus, applying \eqref{important sum for the final estimate} and \eqref{estimate for E-+ delta}, we achieve the following
    \begin{align}
        \ln|\det (P^\delta(q) - z)| -  N^d \left( \int_{\mathbb{T}^{2d}} \ln \, |p_0(\rho) - z| d\rho \! \right) 
        &\leq  N^d \! \left(  \mathcal{O}\left( N^{-\kappa} \ln(N) \right) + \mathcal{O} (N^{-\frac{d}{2}})\right) \!\leq \!N^d  \varepsilon_0(N), \label{inequality for det Pdelta for generic q in the ball}
    \end{align}
    where
    \begin{equation} \label{def of eps_0(h)}
        \varepsilon_0(N) := C'' \left(N^{-\frac{d}{2}} + \left( N^{-d} \ln (N) +  N^{-\kappa} \right) \left( \ln \left(\frac{1}{\tau_0}\right) +  \left (\ln (N) \right ) ^2 \right)  \right).
    \end{equation}
    for a sufficiently large constant $C''>0$. Moreover, from Proposition \ref{Prop for bound for q0}, there exists $ q_0 \in \mathbb{C}^{N^d}$ with $ \| q_0 \| _{2} <  \frac{R}{2}$, 
    such that
    \begin{equation} \label{inequality for det Pdelta for q_0}
        \ln|\det (P^\delta(q_0) - z) | - N^d\left(  \int_{\mathbb{T}^{2d}} \ln \, |p_0(\rho) - z| d\rho \right) 
        \geq  - N^{d} \varepsilon_0(N).
    \end{equation}
    \begin{oss} \label{rem on converg of eps 0}
        If $\tau_0$ is not too small, then $ \varepsilon_0(N) \rightarrow0$, for $N \rightarrow + \infty$. 
        For instance, if we consider $\tau_0 \geq e^{-N^{\frac{\kappa}{2}}}$, then for every $m \geq \kappa$ we have $N^{-m} \ln \left( \frac{1}{\tau_0}\right) \leq N^{ \frac{\kappa}{2} - m}$.
        This implies that for each $0 \leq l < \frac{\kappa}{2} $, the following is true
        \begin{equation} \label{conv of eps0(N)}
            \lim_{N \rightarrow + \infty} N^l \varepsilon_0(N) = 0.
        \end{equation}
    \end{oss}
    
    \subsection{Log-determinant estimates under random perturbations}
    In this section we follow the same approach and the same techniques used by Sj\"ostrand in \cite[Chapter 8]{Sjostrand 2} and we adapt his work to our discrete framework.
    
    We fix $z \in \mathbb{C}$ and $\tau_0 \in ] 0, (C_\tau N^{-1})^{1/2} ]$ as in Section \ref{Chap 3} and we define the function
    \begin{align}
        F : 
        \mathbb{C}^{N^d}\ni q \quad & \longmapsto \det (P^\delta (q) -z) \exp \left( -N^d \int_{\mathbb{T}^{2d}} \ln |p_0(\rho) -z| d \rho \right) \in \mathbb{C}, \label{definition of F}
    \end{align}
    where $P^\delta (q) := P + \delta V_q$, with $V_q = \operatorname{diag}(q(n); \, n= 1, \ldots, N^d) \in \mathbb{C}^{N^d \times N^d}$ and $\delta = \frac{\tau_0 N^{-2d}}{C}$, defined in the same way as in the previous section. We observe that $F$ is a holomorphic function in $q$ and, from \eqref{inequality for det Pdelta for generic q in the ball} and \eqref{inequality for det Pdelta for q_0}, we know that
    \begin{align}
        &\ln | F(q) | \leq \varepsilon_0 (N) N^{d} \quad \text{if } \| q \|_{2} < 3 R,  \label{first inequality for F} \\ 
        & \ln | F(q_0) | \geq -\varepsilon_0(N) N^{d}, \quad \text{where } \| q_0 \|_2 <  \frac{R}{2}, \label{second inequality for F}
    \end{align}
    where $\varepsilon_0(N)$ is defined as in \eqref{def of eps_0(h)} and $R>0$ satisfies \eqref{interval of definition of R}. Let $q_1 \in \mathbb{C}^{N^d}$ satisfy $\|q_1\|_{2} = R$. We follow the same strategy as in Section~8 of~\cite{Sjostrand 2}. 
    More precisely, replacing (8.5) and (8.6) in~\cite{Sjostrand 2} 
    with \eqref{first inequality for F} and \eqref{second inequality for F}, 
    respectively, and recalling that $N = (2\pi h)^{-1}$, 
    we may repeat the argument \textit{mutatis mutandis}. This yields the following.

    \begin{prop} \label{prop with estim for set}
        Let $q_1 \in \mathbb{C}^{N^d}$ such that $\|q_1\|_{2} = R$ and assume that $\epsilon > 0$ is small enough so that $\epsilon \exp(\mathcal{O}(\varepsilon_0(N) N^{d}))\leq 1$. Set $r_0>0$ such that $r_0^2 := 1- \left\| \frac{q_0}{R} \right\|_2^2 + \left| \left \langle \frac{q_0}{R},\frac{q_1}{R}\right \rangle\right|^2 $ Then the following is true:
        \begin{equation} \label{prop regarding measure}
            \lambda\left( \{ r \in [0,r_0] : \|q_0 +r q_1\|_{2} < R, \ |F(q_0 +r q_1)| < \epsilon \} \right)  \leq \varepsilon_0(N)N^{d} \exp \left( \mathcal{O}(1) \!+\! \frac{N^{-d}  \ln (\epsilon)}{\mathcal{O}(1)\varepsilon_0(N)} \right),
        \end{equation}
        where the symbol $\mathcal{O}(1)$ in a denominator indicates a bounded positive quantity.
    \end{prop}
    Taking $\epsilon = \exp \left(- \varepsilon_0(N)N^{d +\alpha}\right)$, for a small constant $\alpha > 0$, the upper bound in \eqref{prop regarding measure} is equal to $\varepsilon_0(N)N^{d} \exp \left( \mathcal{O}(1) - \frac{N^{\alpha}}{\mathcal{O}(1)} \right)$.
    
    We shift to a probabilistic setting by introducing a new measure on $\mathbb{C}^{N^d}$. We consider the probability space given by the ball $B_{N^d}(0,R)$ equipped with a probability measure of the form
    \begin{equation}  \label{def of the probability measure}
        \mathbb{P}(dq) = C(N) e^{\Phi (q)} \lambda (dq),
    \end{equation}
    where $\lambda(dq)$ is the Lebesgue measure on $\mathbb{C}^{N^d}$ and $\Phi \in C^1(B_{N^d}(0,R); \mathbb{R})$ may depend on the semi-classical parameter $N$ and satisfies
    
    \begin{equation} \label{estimate on phi}
        \|\nabla \Phi(q)\|_2 = \mathcal{O}(N^{M_\Phi}), \quad \forall\, q \in B_{N^d}(0,R),
    \end{equation}
    for a fixed constant $M_{\Phi} >0$; $C(N)$ is the normalization constant. Again, we use the same arguments used by Sj\"ostrand in Section 8 of \cite{Sjostrand 2} in order to obtain a probabilistic upper bound for $|F|$, and we obtain the following result.
 
    \begin{prop} \label{prob of F(q) < eps}
        Let $\epsilon >0$ be small enough so that $\epsilon \exp(\mathcal{O}(\varepsilon_0(N) N^{d}))\leq 1$. Then
        \begin{equation} \label{final result prop}
            \mathbb{P}(| F(q)| < \epsilon) \leq   \mathcal{O}(1) N^{M_\Phi+3d}\varepsilon_0(N) \exp \left(  \frac{N^{-d}}{\mathcal{O}(1)\varepsilon_0(N)} \ln (\epsilon) \right).
        \end{equation}
    \end{prop}

    \begin{oss} \label{rem on prop 4.6}
        Taking \( \epsilon = \exp\left( -\varepsilon_0(N)N^{d+\alpha} \right)\) for a small constant $\alpha >0$ as before, result \eqref{final result prop} implies that the inequality
        \begin{equation*}
            \ln |\det (P^\delta (q) -z)|  -N^d \int_{\mathbb{T}^{2d}} \ln |p_0(\rho) -z| d \rho  <  - \varepsilon_0(N)N^{d+\alpha}
        \end{equation*}
        is true with probability $\leq  \mathcal{O}(1) \varepsilon_0(N)N^{M_\Phi + 3d} \exp \left(  \frac{-N^{\alpha}}{\mathcal{O}(1)} \right).$
    \end{oss}

    \section{Eigenvalue counting} \label{Counting of the eigenvalues}
    In the final section, we take the last steps for establishing a Weyl law for randomly perturbed operators. We apply two different techniques: first, we derive probabilistic estimates for the number of eigenvalues of the perturbed operator lying in a given region of the complex plane, under suitable regularity assumptions on its boundary; then, we show the weak convergence of the associated eigenvalue counting measures as $N \rightarrow + \infty$.
    
    \subsection[Zero counting techniques for holomorphic functions]{Zero counting techniques for holomorphic functions and applications to spectral estimates} \label{section 5.1}
    We begin by applying a complex-analytic argument: we reduce the problem to counting the zeros of holomorphic functions subject to explicit exponential growth bounds.  
    The key result we rely on is Theorem 1.2 from \cite{zeros of holo fct}. Before stating it, we introduce a notion of Lipschitz boundary, which is required in its formulation. 
    Let $\Gamma \Subset \mathbb{C}$ be an open set and let $\gamma := \partial \Gamma$ be its boundary. Consider $r : \gamma \rightarrow ]0, \infty[$ a Lipschitz function with a Lipschitz modulus $\leq \frac{1}{2}$.
    
    We assume that the boundary $\gamma$ is Lipschitz with respect to the weight $r$ in the following precise sense: we suppose that there exists a constant $C_0>0$ such that for every $x \in \gamma$, there exist new affine coordinates $\widetilde{y}= (\widetilde{y_1}, \widetilde{y_2})$ of the form $\widetilde{y}= U(y-x),$ for $y \in \mathbb{C} \cong \mathbb{R}^2$, being the old coordinates, where $U = U_x$ is an orthogonal matrix, such that the intersection of $\Gamma$ and the rectangle 
    $ R_x:= \{ y \in \mathbb{C} : |\widetilde{y_1}| < r(x), |\widetilde{y_2}| < C_0 r(x)\}$, 
    corresponds to
    \begin{equation} \label{lip weight 2}
        \{ y \in R_x : |\widetilde{y_1}| < r(x), \widetilde{y_2} > f_x(\widetilde{y_1})\}.
    \end{equation}
    Here $f_x(\widetilde{y_1})$ is a Lipschitz function defined on $[-r(x), r(x)]$ with Lipschitz modulus $\leq C_0$. 
    We notice that the Lipschitz assumption on $r$ remains valid if we decrease the weight $r$. Moreover, it will be convenient to extend the function $r$ to all the complex plane $\mathbb{C}$, by defining
    \begin{equation} \label{extension of r}
        r(x) = \inf_{y \in \gamma} \left( r(y) + \frac{1}{2}|x-y|\right).
    \end{equation}
    The extended function is still Lipschitz with modulus $\leq \frac{1}{2}$.
    
    \begin{teo}[Theorem 1.2 of \cite{zeros of holo fct}] \label{Theorem 1.2 Sjostrand}
        Let $\Gamma \Subset \mathbb{C}$ be a simply connected set with a Lipschitz boundary $\gamma$ with an associated Lipschitz weight $r$ as in \eqref{lip weight 2} and \eqref{extension of r}. We put $\widetilde{\gamma}_{ r} := \bigcup_{x \in \gamma} D(x,  r(x))$. Consider $M$ points $z_j^0 \in \gamma$, with $j \in \mathbb{Z}/ M \mathbb{Z}$, distributed along the boundary in the positively oriented sense so that 
        \begin{equation} \label{conditions for points z_j^0}
            \frac{r(z_j^0)}{4} \leq |z_{j+1}^0 - z_j^0| \leq \frac{r(z_j^0)}{2}.
        \end{equation}
        Let $0< h \leq 1$ and let $\phi $ be a continuous subharmonic function on $\widetilde{\gamma}_{r}$ with a distribution extension to $\Gamma \cup \widetilde{\gamma}_{r}$, denoted by the same symbol.
        Then there exists a constant $C_1 >0$ depending only on the constant $C_0$ and there exist $M$ points $\widetilde{z_j} \in D \left(z_j^0, \frac{r(z_j^0)}{2C_1} \right)$, such that if $u $ is a holomorphic function on $\Gamma \cup \widetilde{\gamma}_{r}$, which satisfies
        \begin{align}
            &h \ln |u(z)| \leq \phi(z), \quad \text{on } \, \widetilde{\gamma}_{r}, \\
            &h \ln |u(\widetilde{z_j})| \geq \phi(\widetilde{z_j}) - \epsilon_j, \quad \text{for } \, j = 1,\dots, M,        
        \end{align}
        where $\epsilon_j \geq 0$, then the number of zeros of $u$ in $\Gamma$ satisfies the following inequality
        \begin{align} \label{formula Theorem 1.2 Sjostrand}
            \left | \#(u^{-1}(0) \cap \Gamma) - \frac{1}{2 \pi h} \mu (\Gamma)  \right| \leq \frac{C_2}{h} \left( \mu (\widetilde{\gamma}_{r}) + \sum_{j=1}^M  \epsilon_j \right)
        \end{align}
        for a constant $C_2 > 0$ which only depends on $C_0$ and on $C_1$. In particular, we have taken $\mu := \Delta \phi \in \mathcal{D}'(\Gamma \cup \widetilde{\gamma}_{r})$, which is a positive measure on $\widetilde{\gamma}_{r}$ so that $\mu(\Gamma)$ and $\mu (\widetilde{\gamma}_{ r})$ are well-defined.
    \end{teo}
    Now we apply Theorem \ref{Theorem 1.2 Sjostrand} to our setting. From this point on, we work under the same assumptions stated at the beginning of Section \ref{Chapter 4} and we add the following ones. Let $\Gamma \Subset \mathbb{C}$ be a relatively compact, simply connected open set, independent of $N$, with uniformly Lipschitz boundary $\gamma = \partial \Gamma$ as defined above. We denote by $r_0 > 0$ its associated weight. Furthermore, assume that there exists $\kappa \in ]0,1]$ such that the volume condition \eqref{bounded volume of the operator} holds uniformly for all
    $ z \in \widetilde{\gamma}_r := \gamma + D(0,r),  0 \leq r \leq r_0$.

    Let $\mathbb{P}$ be a probability measure on the ball $B_{N^d}(0,R)$, defined as in \eqref{def of the probability measure} and \eqref{estimate on phi}, with radius $R>0$ which satisfies \eqref{interval of definition of R}. We recall the definition of $\varepsilon_0(N)$ given in \eqref{def of eps_0(h)} and we take $\alpha = C_0 N^{-1}$.
    Finally, we set $P:= p_N$ and $P^\delta (q) := P +\delta V_q$.
    
    We recall the definition of logarithmic potential and some of its properties.
    \begin{defi} [Logarithmic potential] \label{def log pot}
        Let $\mu$ be a finite Borel measure on $\mathbb{C}$ with compact support. We define the logarithmic potential of $\mu$ as
        \begin{equation}
            U_\mu (z) := - \int_{\mathbb{C}} \ln | z-x | \mu(dx), \quad \text{for } z \in \mathbb{C}.
        \end{equation}
    \end{defi}
    We observe that for every finite Borel measure $\mu$ on $\mathbb{C}$ with compact support, we have $U_\mu \in L^1_{\text{loc}}(\mathbb{C}, \lambda(dz))$, and hence $ U_\mu (z) < + \infty, \ \text{ for a.e. } z \in \mathbb{C}$. Moreover,
    \begin{equation} \label{laplacian of log pot}
        \Delta U_\mu =  - 2 \pi \mu, \quad \text{in } \mathcal{D'}(\mathbb{C}).
    \end{equation}
    
    We now consider $(p_0)_*(\lambda)$, the push-forward of the Lebesgue measure $\lambda(d \rho)$ on $\mathbb{T}^{2d}$ under the principal symbol $p_0$ of $p$. Since $p_0 $ is continuous and $\mathbb{T}^{2d}$ is compact, the image ${\operatorname{Im}(p_0)}$ is compact, and consequently $(p_0)_*(\lambda)$ is a finite Borel measure supported on ${\operatorname{Im}(p_0)}$. 
    Accordingly, we set $\phi(z) := - U_{(p_0)_*(\lambda)}(z)$.
    By \eqref{laplacian of log pot} we know that $\phi \in L^1_{\text{loc}}(\mathbb{C}) \subset \mathcal{D}'(\mathbb{C})$ and that its Laplacian 
    \begin{equation} \label{lapl of phi}
        \Delta \phi = 2 \pi (p_0)_*(\lambda), \quad \text{in } \, \mathcal{D}'(\mathbb{C}).
    \end{equation}
    Hence $\phi$ is subharmonic. 
    We now aim to prove that $\phi$ is continuous in a neighborhood of $\gamma$, specifically in  $\widetilde{\gamma_r}$. 
    To do so, we invoke Theorem 1 from \cite{Continuous_Potentials_and_Linear_Mass_Distributions}, without reproducing its proof.
    \begin{teo}[Theorem 1 of \cite{Continuous_Potentials_and_Linear_Mass_Distributions}] \label{theroem about mass distribution}
        Let $m $ be a finite, positive Borel measure with compact support and let $U_m$ be its logarithmic potential. Then $U_m$ is continuous at a point $z_0 \in \mathbb{C}$ if and only if
        \begin{equation}
            \lim_{r \rightarrow 0} \left \{ \limsup\limits_{z\rightarrow z_0} \int_0 ^r \frac{m_t(z)}{t} dt \right\} = 0,
        \end{equation}
        where $m_t(z) := m(D(z,t))$. 
    \end{teo}
    Now, let $z_0 \in \widetilde{\gamma_r}$ and consider $z \in \operatorname{neigh}(z_0) \subset \widetilde{\gamma_r}$. If we take $0 \leq r \ll 1$, by applying the bound \eqref{bounded volume of the operator} uniform  on $\widetilde{\gamma_r}$, we obtain
    \begin{equation*}
          \int_0 ^r \frac{(p_0)_*(\lambda(d \rho))(D(z,t))}{t} dt  = \int_0 ^r \frac{\lambda(\{ \rho \in \mathbb{T}^{2d} : |p_0(\rho) -z|^2 < t^2\})}{t} dt = \int_0 ^r \frac{\mathcal{O}(t^{2\kappa})}{t} dt \leq \frac{C}{2\kappa} r^{2\kappa},
    \end{equation*}
    for a constant $C>0$. Since this result is true uniformly on a neighborhood of $z_0$, it follows that
    \begin{equation*}
        \lim_{r \rightarrow 0} \left \{ \limsup\limits_{z\rightarrow z_0} \int_0 ^r \frac{(p_0)_*(\lambda(d \rho))(D(z,t))}{t} dt \right\} \leq \lim_{r \rightarrow 0} \frac{C}{2\kappa} r^{2\kappa} = 0.
    \end{equation*}
    Therefore, by Theorem \ref{theroem about mass distribution}, $\phi$ is continuous in $\widetilde{\gamma_r}$. 
    
    Now we define, for $q \in B_{N^d}(0,R)$, $u(z) := \det ( P^\delta(q) - z), \ z \in \mathbb{C}$,
    and we observe that $u$ is a holomorphic function on $\mathbb{C}$. 
    Since $\Gamma$ is relatively compact and its boundary $\gamma$ is locally Lipschitz with respect to the weight $r_0$, there exist $l$ points $x_1, \dots, x_l \in \gamma$ such that $ \gamma \subset \bigcup_{j=1}^{l} A_{x_j}$, where
    $$A_{x_j} := \{ y \in R_{x_j} : \widetilde{y} = U_{x_j}(y-x_j), \ |\widetilde{y_1}| < r_0 , \ \widetilde{y_2} = f_{x_j}(\widetilde{y_1})\}, \qquad j= 1, \dots, l.$$
    Here the rectangles $R_{x_j}$, the Lipschitz functions $f_{x_j}$ and the orthogonal matrices $U_{x_j}$ are defined as in the start of this section.
    We parametrize each patch with the function
    $$\Phi_{j}(t) := x_j + U_{x_j}^{-1}((t,f_{x_j}(t))), \qquad t \in ]-r_0, r_0[.$$
    By Rademacher's theorem $f_{x_j}'(t)$ exists for a.e. $t \in ]-r_0, r_0[$ and its modulus is bounded by the Lipschitz constant $L_j$ of $f_{x_j}$. 
    Let $\mathcal{ H}^1_{|\gamma}$ denote the arc length on $\gamma$. From the orthogonality of $U_{x_j}$, we get
    $$\mathcal{H}^1(\gamma) \leq \sum_{j=1}^{l} \int_{-r_0}^{r_0} \| \Phi_j'(t)\|_{2} dt \leq  \sum_{j=1}^{l} 2 r_0 \sqrt{1 + L_j^2} < + \infty,$$
    and consequently, by the hypotheses of $N$-independence, $\mathcal{H}^1(\gamma) = \mathcal{O}(1)$. 
    Now we take a smaller, possibly $N$-dependent weight $0<r\ll 1$ in the definition of Lipschitz boundary and we choose $N_0$ points $ z_j^0 \in \gamma, \ j = 1,2, \dots, N_0,$ 
    satisfying the conditions \eqref{conditions for points z_j^0}. By definition of arc length, we have $ N_0 = \mathcal{O}(r^{-1})$. 
    Applying Proposition \ref{prob of F(q) < eps}, we know that for every fixed $z \in \widetilde{\gamma_r}$, the inequality
    \begin{equation} \label{to be rewritten}
                     |\det (P^\delta (q) -z)|  \geq \epsilon \exp \left( N^d \int_{\mathbb{T}^{2d}} \ln |p_0(\rho) -z| d \rho \right),  
    \end{equation}
    is true with a probability $\geq 1-  \mathcal{O}(1) N^{M_\Phi+3d}\varepsilon_0(N) \exp \left( - \frac{N^{-d}}{\mathcal{O}(1)\varepsilon_0(N)} \ln \left(\frac{1}{\epsilon}\right) \right),$ provided that $\epsilon > 0 $ is small enough so that $\epsilon \exp(\mathcal{O}(\varepsilon_0(N) N^{d}))\leq 1$. 
    
    We define $\widetilde{\epsilon} := N^{-d} \ln \left( \epsilon^{-1} \right), $ so that $\epsilon = \exp (- \widetilde{\epsilon} N^d)$ and we observe that the condition on $\epsilon$ above holds if $\widetilde{\epsilon} \geq C\varepsilon_0 (N)$,
    for some large constant $C>0$. Thus, we can rephrase \eqref{to be rewritten} as it follows. 
    For every $z \in \widetilde{\gamma_r}$, 
    \begin{equation*}
        \mathbb{P} \left(N^{-d} \ln |u(z)| \geq  \phi(z)  - \widetilde{\epsilon} \right) \geq 1-  \mathcal{O}(1) N^{M_\Phi+3d}\varepsilon_0(N) \exp \left( - \frac{\widetilde{\epsilon}}{C\varepsilon_0(N)}  \right).
    \end{equation*}
    We apply Theorem \ref{Theorem 1.2 Sjostrand} and we take $\widetilde{z_1}, \dots, \widetilde{z_{N_0}} \in \widetilde{\gamma_r}$ as in its statement. Applying what we have just found to each of them and we obtain
    \begin{align*}
        \mathbb{P}\left(\bigcap_{j=1} ^{N_0}  \left( N^{-d} \ln |u(\widetilde{z_j})| \geq \phi(\widetilde{z_j})  - \widetilde{\epsilon}\right)\right) \geq  1 - N_0 \mathcal{O}(1) \varepsilon_0(N) N^{M_\Phi +3d} e^{ \left( - \frac{\widetilde{\epsilon}}{C\varepsilon_0(N)}  \right)}.
    \end{align*}
    Furthermore, we recall that inequality \eqref{inequality for det Pdelta for generic q in the ball} obtained in Section \ref{section about estimate for a general potential} is true for all $q \in B_{N^d}(0,3R)$ and for every $z \in \widetilde{\gamma_r}$.
    Hence, in the probability space $(B_{N^d}(0,R), \mathbb{P})$, the inequality
    \begin{equation*}
        N^{-d} \ln |u(z)| \leq \phi(z) + \varepsilon_0(N) \leq  \phi(z) + \frac{\widetilde{\epsilon}}{C},
    \end{equation*}
    is true with probability equal to $1$ and for all $z \in \widetilde{\gamma_r}$. 
    Finally, we define the function $\widetilde{\phi}(z) := \phi(z) + \frac{\widetilde{\epsilon}}{C}$, which still satisfies the same hypotheses of $\phi$. Then with probability $\geq  1 - \mathcal{O}(1) r^{-1} \varepsilon_0(N)N^{M_\Phi +3d} \exp \left( - \frac{\widetilde{\epsilon}}{C\varepsilon_0(N)}  \right),$ the following two are true 
    \begin{align*}
           & N^{-d} \ln |u(z)| \leq   \widetilde{\phi}(z), \quad \text{for all } z \in \widetilde{ \gamma_r} \\
            &N^{-d} \ln |u(\widetilde{z_j})| \geq \widetilde{\phi}(\widetilde{z_j})  - \widetilde{\epsilon} \left(1 + \frac{1}{C} \right), \quad \text{for } j = 1, \dots, N_0.
    \end{align*}
    Thus, applying Theorem \ref{Theorem 1.2 Sjostrand} we have proved Theorem \ref{teo on quantitative final result}.
    
    \subsection{Eigenvalue distribution via empirical measure convergence}
    In this section we work within the same framework and use the same definitions as in the previous one. We now adopt a more qualitative perspective based on weak convergence of empirical measures, rather than providing explicit error estimates.
    \begin{defi} [Empirical measure]
        Let $P \in \mathbb{C}^{M \times M}$ be a linear operator, with spectrum $\sigma(P)$. We define the empirical measure of its eigenvalues as the probability measure $ \mu := \frac{1}{M} \sum _{\zeta \in \sigma(P)} \delta_\zeta$.
    \end{defi}
    Hence, the empirical measure for the perturbed operator $P^\delta_N = P^\delta(q)$ is $\xi_N = \frac{1}{N^d}  \sum _{\zeta \in \sigma(P_N^\delta)} \delta_\zeta$. 
    As in Section \ref{section 5.1} we consider the following measure
    \begin{equation}
        \xi:= (p_{0})_{*}(\lambda(d\rho)),
    \end{equation}
    which is the push-forward of the Lebesgue measure $\lambda(d\rho)$ on $\mathbb{T}^{2d}$ through the principal symbol of $p$. 
    By Definition \ref{def log pot}, the logarithmic potentials of the measures $\xi_N, \, \xi$ are given respectively by
    \begin{align} \label{log pot of xiN and xi}
        U_{\xi_N}(z) =  - N^{-d}  \ln \left| \det (P^\delta_N - z) \right|, \qquad
        U_{\xi}(z) = - \int_{\mathbb{T}^{2d}} \ln | p_0(\rho) - z | d\rho.  
    \end{align}

    An additional property of the logarithmic potentials that we will need, is the following: let $\{ \mu_n \}_{n \in \mathbb{N}}$ be a sequence of finite Borel measures with compact support, and suppose that the supports of all $\mu_n$ are contained in a fixed compact set. Then, the almost sure convergence of the associated logarithmic potentials $U_{\mu_n}(z) \to U_\mu(z)$, where $\mu$ is a finite Borel measure with compact support, implies the weak convergence of the measures $\mu_n \rightharpoonup \mu$. 
    
    We have considered this result in the setting of random measures, as stated in \cite[Theorem 2.8.3]{Random matr theory}, and adapted it to our framework.

    We denote $\mathcal{M}(\mathbb{C})$ the space of all finite Borel measures on $\mathbb{C}$ with compact support.
    \begin{teo} \label{Theorem on weak convergence}
        Let $K, K' \Subset \mathbb{C}$ be open, relatively compact sets with $\overline{K} \subset K'$, and let $\{ \mu_n \}_{n \in \mathbb{N}}, \, \mu $ be random measures in $\mathcal{M}(\mathbb{C})$, defined on a probability space $(\Omega, \mathcal{F}, \mathbb{P})$. We assume that all the measures $\mu_n$ and $\mu$ are almost surely uniformly bounded in total mass by the same deterministic constant $B_0>0$ and that there exists $n_0 \geq 1$ sufficiently large such that almost surely
        \begin{equation} \label{supports of the meas mu}
            \operatorname{supp } \mu, \, \operatorname{supp } \mu_n \subset K, \quad \text{for } n \geq n_0.
        \end{equation}
        Suppose that for almost every $z \in K'$, the limit
        \begin{equation}   \label{convergence of the log potentials}
            U_{\mu_n}(z) \rightarrow U_\mu(z), \quad n \rightarrow + \infty, 
        \end{equation}
        is true almost surely (respectively in probability). Then,
        \begin{equation} \label{weak conv of the measures}
            \mu_n \rightharpoonup \mu, \quad n \rightarrow + \infty, \quad \text{weakly}
        \end{equation}
        almost surely (respectively in probability).
    \end{teo}
    \begin{proof}
        We begin with the proof in the case of almost sure convergence. First we notice that the assumption that for a.e. $z \in K'$ the limit \eqref{convergence of the log potentials} holds almost surely is equivalent to assume that almost surely \eqref{convergence of the log potentials} holds for a.e. $z \in K'$. Indeed, if we define $E := \{ (z,\omega) \in K' \times \Omega : U_{\mu_n(\omega)}(z) \rightarrow U_{\mu(\omega)}(z)\} \subset K' \times \Omega$, an application of Tonelli theorem to $\mathbb{1}_E$ proves the claim. 
        
        We observe that $\ln | \cdot - x | \in L^2(K')$ uniformly in $x \in K'$ with a constant $C>0$. Then, by Minkowski's inequality, for almost every $\omega \in \Omega$ and for all $n\geq n_0$, we obtain
        \begin{align*}
             \left( \int_{K'} \left| U_{\mu_n(\omega)}(z)\right|^2 \lambda(dz) \right)^{1/2} \leq   \int_{\mathbb{C}} \left( \int_{K'}| \ln |z-x||^2 \lambda(dz) \right)^{1/2} \! \!\mu_n(\omega)(dx) \leq   \int_{\mathbb{C}} C^{\frac{1}{2}} \mu_n(\omega)(dx) \leq \widetilde{C},
        \end{align*}
        where $\widetilde C>0$ depends only on $K'$ and $B_0$, and is independent of both $n$ and $\omega$.
        By applying a similar argument also for $U_{\mu}$, we conclude that for all $n \geq n_0$ and for almost every $\omega \in \Omega$, $U_{\mu(\omega)}, U_{\mu_n(\omega)}  \in L^2(K')$ uniformly. 
        Combining this result with \eqref{supports of the meas mu} and with the equivalence between the statements proved above, we obtain that there exists a measurable set $\Omega' \subset \Omega$ with $\mathbb{P}(\Omega') = 1$, such that for each $\omega \in \Omega'$, the following are true 
        \begin{enumerate}
            \item $U_{\mu_n(\omega)}(z) \rightarrow U_{\mu(\omega)}(z)$ as $n \rightarrow + \infty$, for a.e. $z \in K'$,
            \item $\operatorname{supp} \mu_n (\omega), \, \operatorname{supp} \mu (\omega) \subset K,$ for all $n \geq n_0$,
            \item there exists a deterministic constant $C_{K',\Omega'}>0$, depending only on $K'$ and $B_0$, and independent of both $n$ and $\omega$, such that for every $\omega \in \Omega'$ and for all $n \geq n_0$,
            \begin{equation} \label{L2 bounds}
            \| U_{\mu_n(\omega)} \|_{L^2(K')}, \| U_{\mu (\omega)} \|_{L^2(K')} \leq C_{K', \Omega'}, 
            \end{equation}
        \end{enumerate}
        In order to prove the statement \eqref{weak conv of the measures}, it is sufficient to show that for all $\omega \in \Omega'$ and for any $\varphi \in C_c^{\infty}(K', \mathbb{R})$,
        \begin{equation} \label{weak converg of log pot}
            \langle \mu_n(\omega) , \varphi \rangle_{\mathcal{D'}, \mathcal{D}} \rightarrow \langle \mu(\omega) , \varphi \rangle_{\mathcal{D'}, \mathcal{D}},  \quad\text{for } n \rightarrow +\infty.
        \end{equation}
        Fix $\omega \in \Omega'$, let $M>0$ and define $ g_n^M(z) := \min (| U_{\mu_n(\omega)}(z)- U_{\mu(\omega)}(z)|, M)$,
        for $z \in K'$ and $n \geq n_0$. By the dominated convergence theorem, we know that
        $g_n^M \underset{n \to \infty}{\longrightarrow} 0,  \text{in } L^1(K'),   \text{for any } M> 0$.
        We use Markov's inequality and the uniform $L^2(K')$-bounds from (3) to estimate:
        \begin{align} \label{big ineq for gnM}
            \left\| g_n^M - |U_{\mu_n(\omega)}-U_{\mu(\omega)}| \right\|_{L^1} 
            \leq \left( \!\int_{\{ |U_{\mu_n(\omega)}(z)-U_{\mu(\omega)}(z)| \geq M\}} \lambda(dz) \! \right )^{\frac{1}{2}} \! \! \left \|  U_{\mu_n(\omega)}-U_{\mu(\omega)}  \right \|_{L^2}
            \leq \frac{4 C^2_{K', \Omega'}} {M}.
        \end{align}
        Hence, for each $\omega \in \Omega'$ and for all $M> 0$, we have
        $ \lim_{n \to +\infty}\| U_{\mu_n(\omega)}-U_{\mu(\omega)} \|_{L^1(K')} \leq \frac{4 C^2_{K', \Omega'}} {M}$.
        Thus, taking the limit for $M \rightarrow + \infty$, we obtain that almost surely 
        $U_{\mu_n} \underset{n \to \infty}{\longrightarrow} U_{\mu}, \ \text{in } \mathcal{D}' (K')$ and using \eqref{laplacian of log pot}, \eqref{weak conv of the measures} holds. 
        
        Now we prove the same result but considering the convergence in probability instead of the almost sure convergence. 
        We first assume that for almost every $z \in K'$,
        and for every $\varepsilon > 0$,
        $ \lim_{n \to \infty} \mathbb{P}\left( | U_{\mu_n}(z) - U_{\mu}(z) | > \varepsilon \right) = 0$.
        Properties (2) and (3) still hold for a measurable set $\Omega' \subset \Omega$ with $\mathbb{P}(\Omega') = 1$. We define $g_n^M(z)$ on $ K'$ as before, recalling also its dependence on $\omega \in \Omega$. We fix $M>0$. From the hypothesis, 
        $$  |g_n^M(z)| \leq M \in L^1(\Omega, \mathbb{P}) \ \text{ and } \ g_n^M(z) \underset{n\to + \infty}{\longrightarrow} 0 \text{ in probability for a.e. } z \in K'.$$
        Thus, for a fixed $\varepsilon >0$, for almost every $z \in K'$, we obtain $ \lim_{n \rightarrow + \infty}\int_\Omega |g_n^M(z)| \mathbb{P}(d\omega) \leq \varepsilon$. Since $\varepsilon $ is arbitrary, we conclude that
        $ g_n^M(z) \underset{n \rightarrow+ \infty}{\longrightarrow} 0,  \text{in }  L^1(\Omega, \mathbb{P})$, for almost every $z \in K'$. Hence, by Tonelli theorem and by the dominated convergence theorem,
       \begin{equation*}
           \int_\Omega \left(\int_{K'} |g_n^M(z)| dz\right) \mathbb{P}(d\omega) \ =\int_{K'} \left(\int_\Omega |g_n^M(z)| \mathbb{P}(d\omega) \right) dz  \underset{n \to +\infty}{\longrightarrow} 0.
       \end{equation*}
       Finally, by Markov's inequality, for any $\varepsilon > 0$
       \begin{equation*}
            \mathbb{P}\left(\int_{K'} |g_n^M(z)| dz > \varepsilon \right) \leq \frac{\int_\Omega \left(\int_{K'} |g_n^M(z)| dz\right) \mathbb{P}(d\omega)}{\varepsilon} \underset{n \to +\infty}{\longrightarrow} 0,
        \end{equation*}
        which proves that $ g_n^M \underset{n \to +\infty}{\longrightarrow} 0  \text{ in } L^1(K'),  \text{ in probability}$,
        for any  $M> 0$.
        We observe that also in this case, inequality \eqref{big ineq for gnM} for $\| g_n^M - |U_{\mu_n}-U_{\mu}| \|_{L^1(K')}$ is still valid almost surely.
        Hence, for every $\varepsilon > 0$ and for each $M > 0$, 
        \begin{align*}
           \mathbb{P}(\|U_{\mu_n}-U_{\mu} \|_{L^1(K')} > \varepsilon) \leq \mathbb{P} \left( \frac{4 C^2_{K', \Omega'}} {M} + \| g_n^M \|_{L^1(K')} > \varepsilon \right) 
        \end{align*}
        Taking $M = 8 C^2_{K', \Omega'} \varepsilon^{-1}$, we obtain $ U_{\mu_n} \underset{n \to +\infty}{\longrightarrow} U_{\mu}  \text{ in } L^1(K'),  \text{ in probability}$. 
        Finally, for every $\varphi \in C_c^\infty(K')$ and for every $\varepsilon> 0$, we get
        \begin{align*}
            \mathbb{P}(|\langle U_{\mu_n}-U_{\mu}, \varphi\rangle_{\mathcal{D'}, \mathcal{D}}| > \varepsilon) \leq \mathbb{P}(\|U_{\mu_n}-U_{\mu}\|_{L^1(K')} \sup_{z \in K'} |\varphi(z)| > \varepsilon) \underset{n \to \infty}{\longrightarrow} 0, 
        \end{align*} 
        which means that $U_{\mu_n} \underset{n \to +\infty}{\longrightarrow} U_{\mu} \, \text{ in } \mathcal{D}' (K')$ in probability and hence 
       $\mu_n \underset{n \to +\infty}{\rightharpoonup} \mu \quad \text{in probability}$.
    \end{proof}
    Finally, we apply Theorem \ref{Theorem on weak convergence} to the measures $\xi_N, \xi$ in the semiclassical limit $N \to + \infty$. Since $\xi_N$ are all probability measures, the hypothesis of the uniform bound on the total mass is satisfied. Recalling that the support of $\xi = (p_0)_*(\lambda) $ is the compact set ${\operatorname{Im}(p_0)}$, we have $ \operatorname{supp} \xi \subset \overline{D}\left(0, \| p_0\|_{L^\infty(\mathbb{T}^{2d})}\right)$.
    Moreover, for each $N \geq 1$, $\operatorname{supp} \xi_N = \sigma(P^\delta_N)\subset \overline{D}\left(0, \left\| P^{\delta}_N\right \| \right)$.
    By Proposition \ref{uniform bound for the norm of p_N,alpha}, there exists a constant $C_1>0$ such that
    \begin{equation*}
         \left\| P^{\delta}_N \right \|  \leq  \| P \| + \delta \|V_q\|\leq   C_1  + \delta \|q\|_{\infty}.
    \end{equation*}
    Since the probability measure $\mathbb{P}$ we are considering is defined on the space $B_{N^d}(0,R)$, we obtain $ \operatorname{supp} \xi_N \subset \overline{D}\left(0, C_1 +  \delta R \right)$,
    where from \eqref{interval of definition of R}, $ \delta R = \mathcal{O}(N^{-1/2}) \leq 1,  \text{ for } N \geq 1$.
    Fix $\eta_0>0$ and define
    $$
    R_*:=\max\!\Big\{\|p_0\|_{L^\infty(\mathbb T^{2d})},\, C_1 + 1\Big\},\qquad
    K:=D(0,R_*+\eta_0)\Subset\mathbb C.
    $$
    Then $\operatorname{supp} \xi, \, \operatorname{supp} \xi_N \subset K$, for all $N \gg 1$.
    We consider an open relatively compact set $K' \Subset \mathbb{C}$ such that $\overline{K} \subset K'$. If we show that for a.e. $z \in K'$, the convergence $ U_{\xi_N}(z) \underset{N \to +\infty}{\longrightarrow}  U_{\xi}(z)$, holds in probability, then the hypotheses of Theorem \ref{Theorem on weak convergence} are satisfied. 
    As observed in \cite[p.~745]{Sjostrand 2}, Sard's theorem implies
    that almost every $z\in\mathbb C$ is a regular value of
    $p_0:\mathbb T^{2d}\to\mathbb C\simeq\mathbb R^2$.
    For each such $z$, \eqref{differentials indep}, applied to $p_0$,
    gives constants $C_z>0$ and $0<t_z\leq1$ such that
    \[
        V_z(t)\leq C_z t\leq C_z t^\kappa,
        \qquad 0<t<t_z,
    \]
    since $0<\kappa\leq1$.
    Thus the volume condition required for the determinant estimates
    holds for almost every $z\in K'$.
    We apply these estimates pointwise in $z$, allowing their constants
    to depend on $z$.
    For each fixed regular value $z$, apply the deterministic construction with $C/a_z$ in place of $C$, where
    $0<a_z\leq1$ will be chosen sufficiently small. Let $u_z$ be the resulting vector and set $q_0(z)=a_z u_z$. For this vector, there exists $D_z >0$, independent of $a_z$, such that $\|u_z\|_2\leq D_zN^{(d-1)/2}$; so choosing $a_zD_z < C_2/2$ ensures $\|q_0\|_2<R/2$. Since \(P + \delta V_{q_0} = P + (a_z\delta) V_{u_z} \), the determinant estimates hold at the original coupling $\delta$, with constants depending on $z$.
    
    Fixing a regular value $z \in K'$ of $p_0$, by \eqref{log pot of xiN and xi} we have $ U_\xi(z) - U_{\xi_{N}}(z) = N^{-d} \ln |F(q)| $, where $F$ is defined as in \eqref{definition of F} and depends on the random vector $q$. Hence, for every $\widetilde{\epsilon} >0$, 
    \begin{equation}
        \mathbb{P} \left( \left|U_\xi(z) - U_{\xi_{N}}(z) \right| > \widetilde{\epsilon} \right) = \mathbb{P} \left( N^{-d} \ln |F(q)| > \widetilde{\epsilon} \right) + \mathbb{P} \left( N^{-d} \ln |F(q)| < - \widetilde{\epsilon}  \right). 
    \end{equation}
    Recalling Remark \ref{rem on converg of eps 0}, we suppose $\tau_0$ is not too small; for instance, we assume $ \tau_0 \geq \exp(-N^{\frac{\kappa}{2}})$.
    
    We fix $\widetilde{\epsilon} > 0$. From \eqref{conv of eps0(N)}, there exists $N_1 = N_1(\widetilde{\epsilon}, z) \geq 0 $ such that $\varepsilon_0(N) \leq \widetilde{\epsilon}, \ \text{for all } \ N \geq N_1$. Applying \eqref{inequality for det Pdelta for generic q in the ball}, for every $N \geq N_1 $ we have
    \begin{equation}\label{1 eg for prob conv}
        \mathbb{P} \left( N^{-d} \ln |F(q)| > \widetilde{\epsilon} \right) \leq \mathbb{P} \left( N^{-d} \ln |F(q)| > \varepsilon_0(N) \right) = 0.  
    \end{equation}
    Moreover, let $0< \alpha < \frac{\kappa}{2} $ and define 
    $ \epsilon_1(N):= \exp(- \varepsilon_0(N)N^{d+\alpha})$.
    With this definition, for $N \gg 1$,
    $ \epsilon_1 \exp(\mathcal{O}(\varepsilon_0(N) N^{d})) \leq 1$. Hence, the hypothesis of Proposition \ref{prob of F(q) < eps} is satisfied with $\epsilon_1$ and from Remark \ref{rem on prop 4.6} it follows
    \begin{equation} \label{prob tending to 0}
        \mathbb{P} \left( N^{-d} \ln |F(q)| < - \varepsilon_0(N)N^{\alpha} \right) \leq  \mathcal{O}(1) \varepsilon_0(N)N^{M_\Phi + 3d} \exp \left(  \frac{-N^{\alpha}}{\mathcal{O}(1)} \right).
    \end{equation}
    Since $0< \alpha < \frac{\kappa}{2}$, Remark \ref{rem on converg of eps 0} ensures the existence of $N_2(z) \geq 0$ such that
    $ \varepsilon_0(N) N^\alpha \leq \widetilde{\epsilon}, \ \text{for all }  N \geq N_2$. Recalling the definition \eqref{def of eps_0(h)}, it follows that the right-hand side of \eqref{prob tending to 0} converges to $0$ as $N\rightarrow + \infty$. Therefore, for any fixed $\eta >0$, there exists $N_3 = N_3(\eta,z) \geq 0$ such that this term is bounded by $\eta$ for every $N \geq N_3$. Consequently, for all $N \geq \max(N_2, N_3)$, we obtain 
    \begin{equation} \label{2 eq for prob conv}
         \mathbb{P} \left( N^{-d} \ln |F(q)| < - \widetilde{\epsilon} \right) \leq \mathbb{P} \left( N^{-d} \ln |F(q)| < - \varepsilon_0(N)N^{\alpha} \right) \leq  \eta.
    \end{equation}
    Finally, combining \eqref{1 eg for prob conv} and \eqref{2 eq for prob conv}, we conclude that for every $\widetilde{\epsilon}>0$ and almost every $z \in K'$, the limit $  \mathbb{P} \left(|U_\xi(z) - U_{\xi_{N}}(z)| >  \widetilde{\epsilon} \right) \xrightarrow[N\rightarrow \infty]{} 0$ is true. Thus, Theorem \ref{Theorem on weak convergence} yields the following convergence in probability:
    \begin{equation} \label{weak conv final result 2}
        \frac{1}{N^d}  \sum _{\zeta \in \sigma(P_N^\delta)} \delta_\zeta\rightharpoonup  (p_0)_*(\lambda(d\rho)), \quad \text{as } \, N \rightarrow + \infty,
    \end{equation}
    which proves Theorem \ref{empirical measure convergence}.
    \begin{oss}
        We observe that the empirical measure $\xi_N$ "counts" the eigenvalues of $P^\delta_N$ over measurable subsets of $\mathbb{C}$. More precisely, for any measurable set $\Gamma \subset \mathbb{C}$,
        \begin{equation} \label{counting eigenvalues}
             \xi_N (\Gamma)= \frac{1}{N^d} \sum _{\zeta \in \sigma(P_N^\delta)} \delta_\zeta(\Gamma) = \frac{1}{N^d} \# \{ \zeta \in \sigma(P^\delta_N) \cap \Gamma\}.
        \end{equation}
    \end{oss}
    By applying a probabilistic version of Portemanteau theorem (see \cite[Theorem 11.1.1]{Dudley}), we obtain that for every measurable set $\Gamma \subset \mathbb C$ such that $ \lambda(p_0^{-1}( \partial \Gamma)) = 0$ and for every $\varepsilon > 0$,
    \begin{equation} \label{fin quantitative result 2}
        \mathbb{P} \left( \left |\frac{1}{N^d}\# \{ \zeta \in \sigma(P^\delta_N) \cap \Gamma\} - \int _{p_0^{-1}( \Gamma)} d \rho  \right | > \varepsilon \right) \longrightarrow 0, \qquad \text{as } N \rightarrow + \infty.
    \end{equation}
    \begin{oss}
        The result in \eqref{fin quantitative result 2} is not substantially different from the quantitative estimate established earlier in \eqref{quantitative final result}. While Theorem \ref{teo on quantitative final result} provides an explicit error bound, but only for a specific class of sets $\Gamma$, \eqref{fin quantitative result 2} yields convergence in probability without an explicit rate and it applies to a broader family of measurable sets.
    \end{oss}
    
    \subsection{Application to the Scottish flag operator} \label{application section}
    We consider a practical example of the results we found. In particular, we work in the same framework we introduced in Section \ref{Main result}. We set $d=1$ and we consider the symbol $p$ defined in \eqref{symbol scottish flag} and its Weyl quantization $p_N$ defined in \eqref{Scott flag op}, which is referred as the \textit{Scottish flag} operator. We note that $p(x,\xi)= p_0(x,\xi) \in S(1,1)$, and it satisfies the symmetry assumption $p(x,\xi) = p(x, -\xi)$.   
    Additionally, we let $\tau_0 = C_0 N^{-m}$ for $C_0>0, m\geq \frac{1}{2}$, so that, by definition \eqref{def of eps_0(h)}, we have
    \begin{align} 
    \varepsilon_0(N) =\mathcal{O} \left( N^{-\frac{1}{2}} + N^{-\kappa}\left (\ln \left(N\right) \right ) ^2   \right),  \quad \text{for }N \gg 1, \label{new eps_0(h)}
    \end{align}
    with $\kappa \in ]0,1]$. 
    By Proposition \ref{4.2.3}, we define $\delta = C^{-1}N^{-m-2}$, for a sufficiently large constant $C>0$.
    Moreover, we set $ C_1 N^{\frac{1 - \kappa}{2}} \leq R \leq C_2 N$, for two constants $C_1, C_2 >0$ and we define $ \mathbb{P}(dq) := C(N) e^{\Phi (q)} \lambda (dq)$, with $\Phi (q) := -\frac{1}{2} \| q\|^2_{2}, \ q \in B_N(0,R)$.
    We observe that
    $\| \nabla \Phi(q)\|_{2} < R = \mathcal{O}(N^{m+\frac{1}{2}})$,
    and thus, by following the definitions in Chapter \ref{Chapter 4}, we obtain
    $M_{\Phi} := m + \frac{1}{2}$.
    We observe that $\sigma(p_N) \subset [-1,1]^2$ and we consider $z_0 \in ]-1,1[^2, R_0 > 0$ such that $\Gamma := D(z_0, R_0)$, with $ \overline{\Gamma} \subset (-1,1)^2 $.
    The set $\Gamma$ satisfies the hypotheses introduced in Section \ref{section 5.1}. For each $z \in \partial \Gamma = \partial D(z_0, R_0)$, we consider affine coordinates $(\widetilde{y_1}, \widetilde{y_2})$ translating $z$ to the origin and rotating so that the tangent to the disc in $z$ is horizontal. Thus, if we denote $y = (y_1,y_2)$ the old coordinates, taking $0 < r_0 < R_0$, the local description of $\Gamma$ in the rectangle $\{ y : (\widetilde{y_1}, \widetilde{y_2}) \in ]-r_0,r_0[ ^2 \}$ is given by
    $ \left \{ y : \widetilde{y_1} \in ]-r_0,r_0[, \widetilde{y_2} > f_z(\widetilde{y_1})\right \}$,
    where the function $f_z(\widetilde{y_1}) := R_0- \sqrt{R_0^2-\widetilde{y_1}^2} $ is Lipschitz in $[-r_0,r_0]$ with modulus $ \leq 1$, whenever $ r_0 \leq \frac{R_0}{\sqrt{2}}$.
    Let $ z \in \mathbb{C}$ and $0 \leq t \ll 1$. By making the change of variables $(u,v) = (\cos(2\pi x), \cos(2\pi\xi))$, we obtain
    \begin{align}
        \lambda \left( \{ (x,\xi) \in \mathbb{T}^2 :|p_0(x,\xi) -z|^2 \leq t\} \right) \! = \! \frac{1}{ \pi ^2} \int_{-1}^{1} \!  \int_{-1}^{1} \! \!  \frac{\mathbb{1}_{\{ (u - \Re (z))^2 + (v - \Im(z))^2 \leq t \}}(u,v)}{\sqrt{(1 \! - \!u^2)(1 \!- \!v^2)}} d u d v. \label{application 1}
    \end{align}
    Now we take a smaller, possibly $N$-dependent weight $0<r\ll1$ in the definition of Lipschitz boundary. In particular, we choose $r$ sufficiently small so that there exists $\eta>0$ with
    \begin{equation} \label{hyp on gamma_r tilde}
    \widetilde{\gamma_r}:= \partial D(z_0,R_0)+D(0,r)\ \subset [-1+\eta, 1 - \eta]^2. 
    \end{equation}
    Hence, taking $t$ sufficiently small, for all $z \in \widetilde{\gamma_r}$, we have
    $$
    \{ (u,v) \in [-1,1]^2 : (u - \Re (z))^2 + (v - \Im (z) )^2 \leq t\} \subset [-1+\delta, 1- \delta]^2,
    $$
    for a constant $\delta >0$ dependent on $\eta$ and $t$. 
    Under these assumptions, from \eqref{application 1}, there exists $C_1>0$ such that
    \begin{align*}
        \lambda \left( \{ (x,\xi) \in \mathbb{T}^2 :|p_0(x,\xi) -z|^2 \leq t\} \right) \leq C_1 \lambda (D(z,\sqrt{t})) = \mathcal{O}(t). 
    \end{align*}
    Therefore the volume condition \eqref{bounded volume of the operator} in the hypotheses of Theorem \ref{teo on quantitative final result} holds with $\kappa = 1$ and by \eqref{new eps_0(h)}, $\varepsilon_0(N) =\mathcal{O} \left( N^{-\frac{1}{2}} \right)$, for $N \gg 1$. We can apply the theorem for any $\widetilde{\epsilon} \geq C_3 \, \varepsilon_0(N)$ (with $C_3 >0$ large enough). For instance, choosing $ 1 \gg \widetilde{\epsilon} := \varepsilon_0(N) N^{\frac{1}{6}} \gg \varepsilon_0(N)$, there exists a constant $C_2>0$ such that
    \begin{align*}
        \left | \# (\sigma(p_N + \delta V_q) \cap \Gamma) - N \lambda \, \left({p^{-1}(\Gamma)} \right) \right| 
        \leq  C_2 N \left( 2 \pi\lambda \, \left(p^{-1}(\partial \Gamma + D(0,r)) \right)  + \mathcal{O}\left(\frac{1}{r} \right) \widetilde{\epsilon} \right),
    \end{align*} 
    with probability  $ \geq 1 - \mathcal{O}(r^{-1}) \varepsilon_0(N)N^{m+\frac{7}{2}} \exp\left(-\frac{N^{1/6}}{C}\right)$.

    Moreover, if $r$ is small enough, it can be shown using the same change of coordinates as before that $ \lambda  \left(p^{-1}(\partial D(z_0,R_0) + D(0,r)) \right) \leq \mathcal{O}\left(r \right)$.
    We choose $\varepsilon_0(N)N^{\frac{1}{6}} \ll r \ll 1$; for instance $r = \varepsilon_0(N) N^{\frac{1}{3}} $. Thus, for a new constant $\widetilde{C_2}>0$, we obtain
     \begin{align*}
        \left | \frac{1}{N} \# (\sigma(p_N + \delta V_q) \cap D(z_0,R_0)) -  \lambda \, \left({p^{-1}(D(z_0,R_0))} \right) \right| 
        \leq  \widetilde{C_2} N^{-\frac{1}{6}} ,
    \end{align*} 
    with probability $ \geq 1 - \mathcal{O}(N^{m+\frac{19}{6}}) e^{- \frac{N^{1/6}}{C}}$,  which tends to $1$ as $N \rightarrow + \infty$. Thus, we have proved Theorem \ref{Scottish operator theorem}.

    \bibliographystyle{abbrvurl}
    \bibliography{references}

\end{document}